\documentclass[conference]{IEEEtran}
\IEEEoverridecommandlockouts
\usepackage{cite}
\usepackage{amsmath,amssymb,amsfonts}
\usepackage{graphicx}
\usepackage{xcolor, soul}
\sethlcolor{green}

\usepackage{booktabs} 
\usepackage{caption,algpseudocode,algorithm,setspace}
\usepackage{scalerel,stackengine}
\usepackage{subfigure}
\usepackage{graphicx}
\usepackage{grffile}
\usepackage{tabu}
\usepackage{textcomp}
 \newcommand{\hatem}[1]{#1}

\def\BibTeX{{\rm B\kern-.05em{\sc i\kern-.025em b}\kern-.08em
    T\kern-.1667em\lower.7ex\hbox{E}\kern-.125emX}}
\begin{document}

\title {Parallel Approximation of the Maximum Likelihood Estimation for the Prediction of Large-Scale Geostatistics Simulations}

\author{\IEEEauthorblockN{Sameh Abdulah, Hatem Ltaief, Ying Sun, Marc G. Genton, and David E. Keyes}
\IEEEauthorblockA{\textit{Extreme Computing Research Center Computer, Electrical, and Mathematical Sciences and Engineering Division,} \\
\textit{King Abdullah University of Science Technology,}\\
Thuwal, Saudi Arabia. \\
Sameh.Abdulah@kaust.edu.sa, Hatem.Ltaief@kaust.edu.sa, Ying.Sun@kaust.edu.sa, \\Marc.Genton@kaust.edu.sa, David.Keyes@kaust.edu.sa.}

}

\maketitle

\begin{abstract}
Maximum likelihood estimation is an important statistical technique for 
estimating missing data, for example in climate and environmental applications, 
which are usually large and feature data points that are irregularly spaced. 
In particular, the Gaussian log-likelihood function
is the \emph{de facto} model, which operates on the resulting sizable dense covariance matrix. The advent of high performance 
systems with advanced computing power and memory capacity 
have enabled full simulations only for rather small dimensional climate problems, solved at the 
machine precision accuracy. The challenge for high dimensional problems lies in the computation requirements 
of the log-likelihood function, which necessitates ${\mathcal O}(n^2)$ storage and
${\mathcal O}(n^3)$ operations, where $n$ represents the number of given spatial locations.
This prohibitive computational cost may be reduced by using 
approximation techniques that not only enable large-scale simulations otherwise intractable,
but also maintain the accuracy and the fidelity of the spatial statistics model.
In this paper, we extend the Exascale GeoStatistics software framework 
(i.e., \texttt{ExaGeoStat}\footnote{https://github.com/ecrc/exageostat}) 
to support the Tile Low-Rank (TLR) approximation technique, which exploits 
the data sparsity of the dense covariance matrix 
by compressing the off-diagonal tiles up to a user-defined accuracy threshold. The 
underlying linear algebra operations may then be carried out on this data compression
format, which may ultimately reduce the arithmetic complexity of the maximum
likelihood estimation and the corresponding memory footprint. Performance results of TLR-based
computations on shared and distributed-memory systems attain up to $13$X and $5$X speedups, 
respectively, compared to full accuracy simulations using synthetic and real datasets (up to $2$M), while 
ensuring adequate prediction accuracy. 
\end{abstract}

\begin{IEEEkeywords}
massively parallel algorithms, machine learning algorithms, applied computing mathematics and 
statistics, maximum likelihood optimization, geo-statistics applications
\end{IEEEkeywords}

\section{Introduction}
Current massively parallel systems provide  unprecedented computing power
with up to millions of execution threads. This hardware technology evolution 
comes at the expense of a limited memory capacity per core, which may prevent 
simulations of big data problems. In particular, climate/weather simulations
usually rely on a complex set of Partial Differential Equations (PDEs) to 
estimate conditions at specific output points based on
semi-empirical models and assimilated measurements.
This conventional approach translates the original big data problem into a 
large-scale simulation problem, solved globally, en route to particular quantities 
of interest, and it relies on PDE solvers to extract performance from the targeted 
architectures. 

An alternative available in many use cases is to estimate missing 
quantities of interest from a statistical model. Until recently, the computation 
used in statistical models, like using field data to estimate parameters of a 
Gaussian log-likelihood function and then evaluating that distribution to 
estimate phenomena where field data are not available, was intractable for very 
large meteorological and environmental datasets. This is due to the arithmetic complexity, 
for which a key step grows as the cube of the problem size~\cite{michael1999interpolation }, 
i.e., increasing the problem 
size by a factor of $10$ requires $1,000$X more work (and $100$X more memory). 
Therefore, the existing hardware landscape with its limited memory capacity, and even
with its high thread concurrency, still appears unfriendly for large-scale
simulations due to the aforementioned curse of dimensionality~\cite{abdulah2017exageostat}.
Piggybacking on the renaissance in hierarchically low rank computational linear algebra,
we propose to exploit data sparsity in 
the resulting, apparently dense, covariance matrix by compressing the off-diagonal blocks 
up to a specific application-dependent accuracy.

This work leverages our Exascale GeoStatistics software framework
(\texttt{ExaGeoStat}~\cite{abdulah2017exageostat}) in the context of climate and environmental simulations, 
which calculates the core statistical operation, i.e.,
the Maximum Likelihood Estimation (MLE),
up to the machine precision accuracy for only rather small spatial datasets. \texttt{ExaGeoStat} relies on
the asynchronous task-based dense linear algebra library \texttt{Chameleon}~\cite{chameleon-soft}
associated with the dynamic runtime system \texttt{Star\-PU}~\cite{augonnet2011starpu}
to exploit the underlying computing power \hatem{toward} large-scale systems.

However, herein, we reduce the memory footprint 
and the arithmetic complexity of the MLE to alleviate the dimensionality bottleneck.
We employ the Tile Low-Rank (TLR) data format for the compression, as implemented
in the Hierarchical Computations on Manycore Architectures 
(\texttt{HiCMA}\footnote{https://github.com/ecrc/hicma})
numerical library. \texttt{HiCMA} relies on task-based
programming model and is deployed on shared~\cite{akbudak-isc17} and 
distributed-memory systems~\cite{akbudak-europar18} via \texttt{StarPU}.
The asynchronous execution achieved via \texttt{StarPU} is even more critical for 
\texttt{HiCMA}'s workload, characterized by lower arithmetic intensity, since it permits to mitigate
the latency overhead engendered by the data movement. 




 \subsection{Contributions}

The contributions of this paper are sixfold.

\begin{itemize}
\item We propose an accurate and amenable MLE
framework using TLR-based approximation format to reduce the prohibitive complexity
of the apparently dense covariance matrix computation.

\item We provide a TLR solution for the prediction operation to impute values related
 to non-sampled locations.

\item We demonstrate the applicability of our approximation technique on both synthetic 
(up to 1M locations) and real datasets (i.e., soil moisture from the Mississippi Basin 
area and wind speed from the Middle East area).

\item We port the \texttt{ExaGeoStat} simulation framework to a myriad of shared and distributed-memory 
systems using a single source code to enhance user productivity, 
thanks to its modular software stack.

\item We conduct a comprehensive performance evaluation to highlight the effectiveness 
of the TLR-based approximation method compared to the original full accuracy approach.
The experimental platforms include shared-memory Intel Xeon Haswell / Broadwell / KNL / Skylake 
high-end HPC servers and the distributed-memory Cray XC40 \emph{Shaheen-2} supercomputer. 

\item We perform a thorough qualitative analysis to 
assess the accuracy of the estimation of the Mat\'{e}rn covariance parameters
as well as the prediction operation. 
Performance results of TLR-based 
MLE computations on shared and distributed-memory systems achieve up to $13$X and $5$X speedups, 
respectively, compared to full machine precision accuracy using synthetic and real environmental datasets (up to $2$M), without 
compromising the prediction quality. 
\end{itemize}
\hatem{The previous works}~\cite{akbudak-isc17,akbudak-europar18} \hatem{focus solely on the standalone 
linear algebra operation, i.e., the Cholesky factorization. They assess its performance
using a simplified version of the Mat\'{e}rn kernel on synthetic datasets. Herein, 
we integrate and leverage these works into the parallel approximation of the 
maximum likelihood estimation for the prediction of large-scale geostatistics simulations.} 



\hatem{The remainder of the paper} is organized as follows.
Section~\ref{sec:related} covers different MLE approximation techniques that have been proposed in the literature.
Section~\ref{sec:motivation} illustrates the climate modeling structure used as a backbone for this work.
Section~\ref{sec:matern_kernel} recalls the necessary background on the Mat\'{e}rn covariance functions. 
Section~\ref{sec:tlr} describes the Tile Low-Rank (TLR) approximation technique and the \texttt{HiCMA} TLR approximation library and its integration into the \texttt{ExaGeoStat} framework.
Section~\ref{sec:exageostat} highlights
the \texttt{ExaGeoStat} framework with its software stack.
Section~\ref{sec:geospat} defines both the synthetic datasets and  the two climate datasets obtained
from large geographic regions, i.e., the Mississippi River Basin and the 
Middle-East region that are used to evaluate the proposed TLR method.
Performance results and accuracy analysis are presented in Section~\ref{sec:perf}, using both synthetic 
and real environmental datasets, and we conclude in Section~\ref{sec:summary}.


\section{Related Work}
\label{sec:related}
Approximation techniques to reduce arithmetic complexities and memory 
footprint for large-scale climate and environmental applications are
well-established in the literature.
Sun et al.~\cite{sun2012} have discussed several of these methods 
such as Kalman filtering~\cite{wikle1999dimension}, moving averages~\cite{ver2004flexible}, 
Gaussian predictive processes~\cite{banerjee2008gaussian}, fixed-rank
kriging~\cite{cressie2008fixed}, covariance tapering~\cite{Kaufman2008,sang2012full}, and 
low-rank splines~\cite{gu2013smoothing}. All these methods depend on low-rank 
models, where a latent process is used with lower dimension, and eventually result
in a low-rank representation of the covariance matrix.
Although these
methods propose several possibilities to reduce the complexity of generating
and computing the domain covariance matrix, several restrictions limit their
functionality~\cite{Stein2013, huang2017hierarchical}. 

On the other hand, low-rank off-diagonal matrix approximation techniques 
have gained a lot of attention to cope with covariance matrices of high dimension.
In the literature, these are commonly
referred as hierarchical matrices or 
$\mathcal{H}$-matrices~\cite{hackbusch1999sparse,khoromskij2009application}. 
The development of various data compression techniques such
as Hierarchically Semi-Separable (HSS)~\cite{ghysels2016efficient}, 
$\mathcal{H}^2$-matrices~\cite{pouransari2015fast,borm2015approximation,sushnikova2016compress},
Hierarchically Off-Diagonal Low-Rank (HODLR)~\cite{aminfar2016fast},
Block/Tile Low-Rank (BLR/TLR)~\cite{amestoy2017complexity,pichon2017sparse,akbudak-isc17,akbudak-europar18}
increases their impact on a wide range of scientific
applications. 
%
Each of the aforementioned data compression formats has pros and cons in terms
of arithmetic complexity, memory footprint and efficient parallel implementation, 
depending on the application operator. We have chosen to rely on TLR 
data compression format, as implemented in the \texttt{HiCMA} library. TLR
may not be the best in terms of theoretical bounds for asymptotic sizes. However, 
thanks to its flat data structure to store the low-rank off-diagonal matrix, TLR
is more versatile to run on various parallel systems. This may not be the case
for data compression formats \hatem{(i.e., $\mathcal{H}$/$\mathcal{H}^2$-matrices, HSS, and HODLR)} 
with recursive tree structure based on nested and non-nested
bases, especially when targeting distributed-memory systems.
\hatem{It is also noteworthy to mention the differences between BLR and TLR. While
these data formats are conceptually identical, BLR has been developed in the context
of multifrontal sparse direct solvers (i.e., MUMPS}~\cite{mumps}\hatem{). The MUMPS-BLR variant takes only dense input matrices 
(i.e., the fronts), computes the Schur complement, and compresses on-the-fly individual blocks
once all their updates have been applied. Therefore, the MUMPS-BLR variant reduces the algorithmic complexity but
not the memory footprint, which may not be a problem since the size of the fronts are typically
much smaller than the problem size. In contrast, the TLR variant in \texttt{HiCMA} accepts dense or already compressed
matrices as inputs, and therefore, permits to reduce arithmetic complexity and memory footprint. The latter is paramount
when operating on dense matrices.}






\section{MLE-Based Climate Modeling and Prediction}
\label{sec:motivation}
Climate and environmental datasets consist of a set of locations regularly
 or irregularly distributed across
a specific geographical region where each location is
associated with a single read of a certain climate and environmental 
variable, for example, wind speed, air pressure, soil moisture, and humidity.


They are usually modeled in geostatistics as a realization from a 
Gaussian spatial random field. Specifically, let ${\bf Z}=\{Z({\bf s}_1),\ldots,Z({\bf s}_n)\}^\top$ 
 be a realization of a Gaussian random field $Z({\bf s})$ at a set of $n$ spatial locations 
${\bf s}_1,\ldots,{\bf s}_n$ in $\Bbb{R}^d$, 
$d\geq 1$. We assume the mean of the random field $Z({\bf s})$  is zero for simplicity and the 
stationary covariance function has a parametric from 
$C({\bf h};{\boldsymbol \theta})=\mbox{cov}\{Z({\bf s}),Z({\bf s}+{\bf h})\}$,
where ${\bf h}\in\Bbb{R}^d$ is a spatial lag vector and ${\boldsymbol \theta}\in\Bbb{R}^q$ is
an unknown parameter vector of interest. 
Denote by ${\boldsymbol \Sigma}({\boldsymbol \theta})$ the covariance matrix with
entries ${\boldsymbol \Sigma_{ij}}=C({\bf s}_i-{\bf s}_j;{\boldsymbol \theta})$, $i,j=1,\ldots,n$.
The matrix ${\boldsymbol \Sigma}({\boldsymbol \theta})$ is symmetric and positive definite.
Statistical inference about $\boldsymbol \theta$ is often based on the Gaussian
log-likelihood function as follows:
\begin{equation}
	\label{eq:likeli}
	\ell({\boldsymbol \theta})=-\frac{n}{2}\log(2\pi) - \frac{1}{2}\log |{{\boldsymbol \Sigma}({\boldsymbol \theta})}|-\frac{1}{2}{\bf Z}^\top {\boldsymbol \Sigma}({\boldsymbol \theta})^{-1}{\bf Z}.
\end{equation}
The main goal is to compute $\widehat{\boldsymbol \theta}$, which represents the maximum
likelihood estimator of ${\boldsymbol \theta}$ in equation (\ref{eq:likeli}). In the case of
large-scale applications, i.e., $n$ is large and locations are irregularly distributed across the region, 
the evaluation of equation (\ref{eq:likeli})
becomes computationally challenging. The log determinant and linear solver involving an $n$-by-$n$ dense 
and unstructured covariance matrix
${\boldsymbol \Sigma}({\boldsymbol \theta})$ require ${\mathcal O}(n^3)$ floating-point operations (flops)
on ${\mathcal O}(n^2)$ memory. Herein lies the challenge. For example, assuming a dataset
  on a grid with approximately $10^3$
 longitude values and $10^3$ latitude values, 
the total number of locations will be $10^6$. Using double-precision floating-point arithmetic, 
the total number of flops will be then equal to one Exaflop with a corresponding memory footprint of 
$10^{12} \times 8 $ bytes $ \sim 80 $ TB, which renders the simulation impossible.


Once $\widehat{\boldsymbol \theta}$ has been computed, 
we can use it to predict unknown measurements 
at a given set of new locations (i.e., supervised learning).
For instance, we can predict $m$ unknown measurements $\mathbf{Z}_1$, where $\mathbf{Z}_2$ represents
a set of $n$ known measurements instead. Thus, the problem can be represented as a multivariate normal
joint distribution~\cite{cressie2015statistics, genton2007separable} as follows
\begin{equation}
	\begin{bmatrix}
		\mathbf{Z}_1 \\
		\mathbf{Z}_2 \\
	\end{bmatrix}
	\sim
	N_{m+n}
\left(
	\begin{bmatrix}
		\boldsymbol{\mu}_{1} \\
		\boldsymbol{\mu}_{2} \\
	\end{bmatrix}
	,
	\begin{bmatrix}
		\mathbf{\Sigma}_{11} & \mathbf{\Sigma}_{12} \\
		\mathbf{\Sigma}_{21} & \mathbf{\Sigma}_{22} \\
	\end{bmatrix}
	\right),
	\label{eq:pred1}
\end{equation}
with $\mathbf{\Sigma}_{11} \in \mathbb{R}^{m \times m}$, $\mathbf{\Sigma}_{12} \in \mathbb{R}^{m \times n}$, $\mathbf{\Sigma}_{21} \in \mathbb{R}^{n \times m}$, and $\mathbf{\Sigma}_{22} \in \mathbb{R}^{n \times n}$. 
The associated conditional distribution can be represented as
\begin{equation}
	\label{eq:pred2}
	\mathbf{Z}_{1}|\mathbf{Z}_{2} \sim N_{m} 
	(
	\boldsymbol{\mu}_{1}+\mathbf{\Sigma}_{12}\mathbf{\Sigma}_{22}^{-1} (\mathbf{Z}_{2}-\boldsymbol{\mu}_{2})
	,
	\mathbf{\Sigma}_{11}-\mathbf{\Sigma}_{12}\mathbf{\Sigma}_{22}^{-1}
	\mathbf{\Sigma}_{21}
	).
\end{equation}
Assuming that the observed vector $\mathbf{Z}_2$ has a zero-mean function 
(i.e., $\boldsymbol{\mu}_{1} = \boldsymbol{0}$ and $\boldsymbol{\mu}_{2} = \boldsymbol{0})$,  
the unknown vector $\mathbf{Z}_1$ can be predicted~\cite{genton2007separable} by solving
\begin{equation}
	\mathbf{Z}_{1}= \mathbf{\Sigma}_{12} \mathbf{\Sigma}_{22}^{-1} \mathbf{Z}_{2}.
	\label{eq:pred3}
\end{equation}
Equation (\ref{eq:pred3}) also depends on two covariance 
matrices, i.e., $\mathbf{\Sigma}_{12}$ and  $\mathbf{\Sigma}_{22}$.
Thus, the prediction operation is as challenging as the initial $\widehat{\boldsymbol \theta}$ estimation operation,
since it also necessitates the Cholesky factorization, followed by a forward and backward
substitution applied on several right-hand sides. 

In this study, we aim at exploiting the data sparsity of the various covariance
matrices by applying a TLR-based approximation technique to reduce the
arithmetic complexity and memory footprint of both operations, i.e.,  the MLE and the prediction,
in the context of climate and environmental applications. 


\section{Mat\'{e}rn Covariance Function}
\label{sec:matern_kernel}

To construct the covariance matrix ${\boldsymbol \Sigma}({\boldsymbol \theta})$ in equation
(\ref{eq:likeli}), a valid (positive definite) parametric covariance model is needed. 
Among the many possible covariance models in the literature, the Mat\'{e}rn 
family~\cite{Matern1986a} has proved useful in a wide 
range of applications. The class of Mat\'{e}rn covariance functions
 is widely used in geostatistics and spatial statistics ~\cite{chiles2009geostatistics}, machine learning~\cite{BoermGarcke2007}, 
image analysis, weather forecasting and climate science. 
Handcock and Stein~\cite{Handcock1993a} introduced the Mat\'{e}rn form of 
spatial correlations into statistics as a flexible parametric class where one parameter determines
the smoothness of the underlying spatial random field. The history of this family
 of models can be found in~\cite{Guttorp2006a}.
The Mat\'{e}rn form also naturally describes the correlation among temperature 
fields that can be explained by simple energy balance climate models~\cite{North2011a}. 
The Mat\'{e}rn class of covariance functions is defined as
\begin{equation}
	\label{eq:MaternCov}
	C(r;{\boldsymbol \theta})=\frac{\theta_1}{2^{\theta_3-1}\Gamma(\theta_3)}\left(\frac{r}{\theta_2}\right)^{\theta_3} {\mathcal K}_{\theta_3}\left(\frac{r}{\theta_2}\right),
\end{equation}
where $r=\|{\bf s}-{\bf s}'\|$ is the distance between two spatial locations, ${\bf s}$ and ${\bf s}'$, and 
${\boldsymbol \theta}=(\theta_1,\theta_2,\theta_3)^\top$.  Here $\theta_1>0$ is the variance, 
$\theta_2>0$ is a spatial range parameter that measures how quickly the correlation of the random 
field decays with distance, and $\theta_3>0$ controls the smoothness of the 
random field, with larger values of $\theta_3$ corresponding to smoother fields. 
the spatial range $\theta_2$ parameter usually can be 
represented by using $0.03$ for weak correlation, $0.1$ for medium correlation, 
and $0.3$ for strong correlation. The smoothness $\theta_3$ parameter,
which represents the data smoothness can be represented by $0.5$ for a rough process,
and $1$ for a smooth process~\cite{genton2017hierarchical}.

The distance between any two spatial locations can be efficiently computed using Euclidian distance. 
However, in the case of real datasets on the surface of a sphere, the Great-Circle Distance (GCD) metric is more suitable. 
The best representation of the GCD
distance is the haversine formula given in~\cite{rick1999deriving}
 \begin{equation}
	\label{eq:gcd}    
{\displaystyle \operatorname {hav} \left({\frac {d}{r}}\right)=\operatorname {hav} (\varphi _{2}-\varphi _{1})+\cos(\varphi _{1})\cos(\varphi _{2})\operatorname {hav} (\lambda _{2}-\lambda _{1})},
\end{equation}
where hav is the haversine function
$\operatorname {hav} (\mathbf \theta )=\sin ^{2}\left({\frac {\mathbf \theta }{2}}\right)={\frac {1-\cos(\mathbf \theta )}{2}}$, 
$d$ is the distance between the two locations, $r$ is the radius of the sphere, 
$\varphi_1$ and $\varphi_2$ are the latitude of location 1 and latitude of 
location 2, in radians, respectively, and $\lambda _{1}$ and $\lambda _{2}$ 
are the counterparts for the longitude.

The function ${\mathcal K}_{\theta_3}$ denotes the modified Bessel function of 
the second kind of order $\theta_3$. When $\theta_3=1/2$, the Mat\'{e}rn 
covariance function reduces to the exponential covariance model 
$C(r;{\boldsymbol \theta})=\theta_1 \exp(-r/\theta_2)$, and describes a 
rough field, whereas when $\theta_3=1$, the Mat\'{e}rn covariance function reduces to
 the Whittle covariance model 
$C(r;{\boldsymbol \theta})=\theta_1 (r/\theta_2){\mathcal K}_1(r/\theta_2)$, 
and describes a smooth field.
The value $\theta_3=\infty$ corresponds to a Gaussian covariance model,
which describes a very smooth field infinitely mean-square differentiable. 
Realizations from a random field with Mat\'{e}rn covariance functions are 
$\lfloor \theta_3-1 \rfloor$ times mean-square differentiable. Thus, the parameter $\theta_3$ 
is used to control the degree of smoothness of the random field. 

In theory, the three parameters of the Mat\'{e}rn covariance function need to be positive real
numbers. Empirical values derived from the empirical covariance of the data can serve
as starting values and provide bounds for the optimization. Moreover, the parameter
$\theta_3$ is rarely found to be larger than 1 or 2 in geophysical applications, as those already
 correspond to very smooth realizations.

\section{Tile Low-Rank Approximation}
\label{sec:tlr}

Tile algorithms have been used for the last decade on manycore architectures
to speedup parallel linear solvers algorithms, as implemented in 
the \texttt{PLASMA} library~\cite{Agullo_2009_jpcs}. 
Compared to \texttt{LAPACK} block algorithms~\cite{anderson1999lapack},
tile algorithms permit to bring the parallelism within multi-threaded BLAS 
to the fore by splitting the matrix into dense tiles. The resulting fine-grained 
computations weaken the synchronizations points and create opportunities for look-ahead
to maximize the hardware occupancy.
In this study, we propose an MLE optimization framework, which operates on 
Tile Low-Rank (TLR) data compression format, as implemented in the
Hierarchical Computations on Manycore Architectures (\texttt{HiCMA}) 
library. \hatem{More details about algorithmic complexity and memory footprint
can be found in}~\cite{akbudak-isc17,akbudak-europar18}.

 \begin{figure}
\centering
\includegraphics[width=0.5\linewidth]{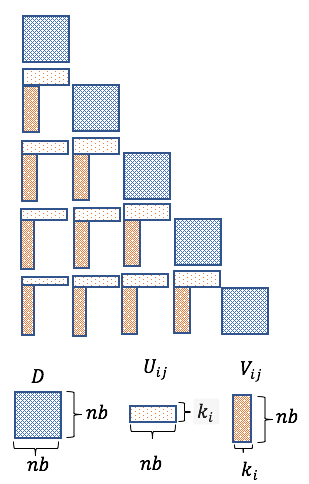}
  \caption{TLR representation of a covariance matrix ${\boldsymbol \Sigma}({\boldsymbol \theta})$ with fixed accuracy.}
  \label{fig:lr_representation}
\end{figure}

 Figure~\ref{fig:lr_representation} illustrates the TLR representation of a given covariance matrix ${\boldsymbol \Sigma}({\boldsymbol \theta})$. Following the same principle as dense tile algorithms~\cite{abdulah2017exageostat}, our covariance 
matrix ${\boldsymbol \Sigma}({\boldsymbol \theta})$ is divided into a set of square tiles. 
\hatem{The Singular Value Decomposition (SVD), Randomized SVD (RSVD), or Adaptive Cross 
Approximation (ACA) may be used then} to approximate each off-diagonal 
tile up to a user-defined accuracy threshold. This threshold is, in fact, application-dependent
and enables to truncate and keep the most significant $k$ singular values and their associated left and right singular
vectors, ${\boldsymbol U}$ and ${\boldsymbol V}$, respectively. The number $k$ is the actual rank and is determined on a tile basis, i.e., 
one should expect variable ranks across the matrix tiles.
A low accuracy translates into small ranks (i.e., low memory footprint), and therefore, brings the arithmetic intensity 
of the overall algorithm close to the memory-bound regime. Conversely, a high accuracy generates
large ranks (high memory-footprint), which increases the arithmetic intensity and makes the algorithm 
run in the compute-bound regime.  
Each tile $(i,j)$ can then be represented by the product of $U_{ij}$ and $V_{ij}$,
with a size of $nb \times k$, where $nb$ represents the tile size. The tile size is a tunable 
parameter and has a direct impact on the overall performance, since it corresponds to the trade-off
between arithmetic intensity and degree of parallelism.

The next section introduces the new extension of {\texttt ExaGeoStat} framework~\cite{abdulah2017exageostat}
toward TLR matrix approximations of the Mat\'{e}rn covariance functions and TLR matrix computations
using the high performance \texttt{HiCMA} numerical library, in the context of MLE calculations.

\section{\texttt{ExaGeoStat} Software Infrastructure}
\label{sec:exageostat}
This work is an extension of our {\texttt ExaGeoStat} framework~\cite{abdulah2017exageostat},
a high performance framework for geospatial statistics in climate and
environment modeling. In~\cite{abdulah2017exageostat}, 
we propose using full machine precision accuracy for maximum likelihood estimation.
 Besides demonstrating the
hardware portability of the framework, one of the motivations
was to provide a reference implementation for eventual performance and accuracy 
assessment against different approximation techniques.
In this work, we extend the {\texttt ExaGeoStat} framework with a TLR approximation
technique and assessing it with the full accuracy reference solution.
{\texttt ExaGeoStat} sits on top of three main components: (1) the {\texttt ExaGeoStat} operational
routines which generate synthetic datasets (if needed), solve the MLE problem,
and predict missing values at non-sampled locations, (2) linear algebra libraries, i.e.,
\texttt{HiCMA} and \texttt{Chameleon},  
which provide linear solvers support for the {\texttt ExaGeoStat} routines, and (3) 
a dynamic runtime system, i.e., \texttt{StarPU}~\cite{augonnet2011starpu}, which orchestrates the adaptive execution
of the main computational routines on various shared and distributed-memory systems.


Once the {\texttt ExaGeoStat} high-level tasks (i.e., matrix generation, log-determinant 
and solve operations) are defined with their respective data dependencies, \texttt{StarPU}
can enroll the sequential code and may asynchronously schedule the various tasks
on the underlying hardware resources. By abstracting the hardware complexity via the 
usage of \texttt{StarPU}, user-productivity may be enhanced since developers can focus
on the numerical correctness of their sequential algorithms and leave the parallel 
execution to the runtime system.
Thanks to an out-of-order execution, \texttt{StarPU} is also capable of reducing
load imbalance, mitigating data movement overhead, and increasing occupancy
on the hardware.

Last but not least, {\texttt ExaGeoStat} includes an optimization software layer
(i.e., The NLopt library) to optimize the likelihood objective function. 

%
%
\section{Definitions of Synthetic and Real Datasets}
\label{sec:geospat}
In this study, we use both synthetic and real datasets to validate our proposed TLR {\texttt ExaGeoStat} 
on different hardware architectures.

 Synthetic data are generated at irregular locations over a predefined region  over a two-dimensional space 
 ~\cite{huang2017hierarchical,sun2016statistically}. The synthetic representation aims at generating 
 spatial locations where no two locations are too close. The data 
 locations are generated using $n^{1/2}(r-0.5+X_{rl},l-0.5+Y_{rl})$ for 
$r,l \in \{1,\dots, n^{1/2}\}$, where $n$ represents the number of locations, and $X_{rl}$ and $Y_{rl}$ 
are generated using a uniform distribution on ($-$0.4, 0.4).

\begin{figure}
\centering
\includegraphics[width=0.6\linewidth ]{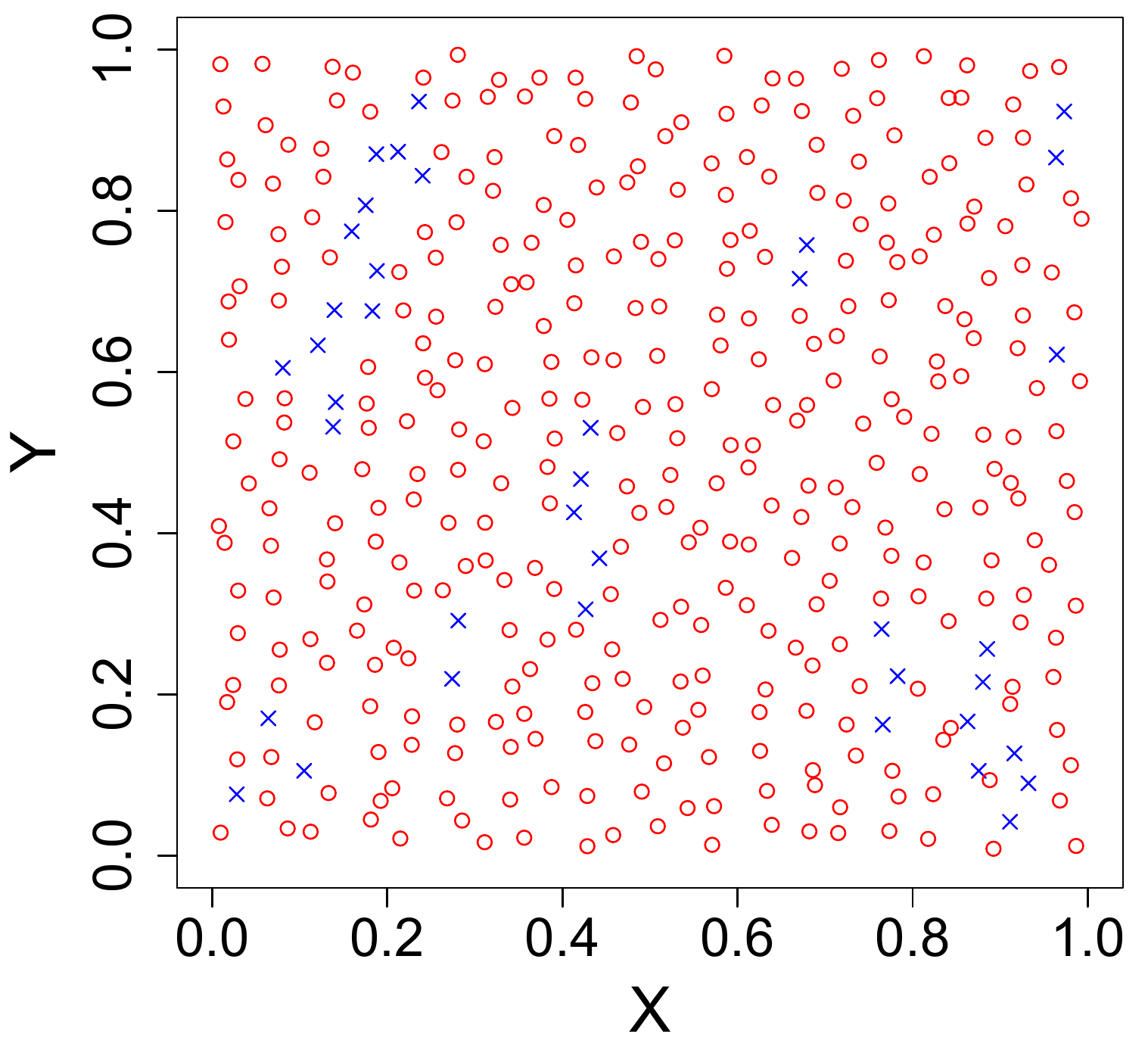}
	\caption{An example of $400$ points irregularly distributed in space, with 362 points (\textcolor{red}{$\circ$}) 
	for maximum likelihood estimation and $38$ points (\textcolor{blue}{$\times$}) for prediction validation.}
	\label{fig:irregular-space}
\end{figure}

A drawable example of 400 irregularly spaced grid locations in a square region is shown by 
Figure~\ref{fig:irregular-space}.  We only use such a small example to highlight how we are
 generating spatial locations, however,  this work uses synthetic datasets up to $ 10^{6}$ locations
  with a total covariance matrix size equals $10^{12}$ double precision elements -- about $80$ TB of memory.

Numerical models are essential tools to improve the understanding of the
global climate system and the cause of different climate variations. Such numerical models
are able to well describe the evolution of many variables related to the climate system, for instance,
temperature, precipitation, humidity, soil moisture, wind speed and pressure, through solving
a set of equations.  The process involves physical parameterization, initial condition configuration, 
numerical integration, and data output.  In this study, we use the proposed methodology
 to investigate the spatial variability of two different
 kinds of climate variables: soil moisture and wind speed.

Soil moisture is a key factor in evaluating the state of the hydrological process and has a wide range
of applications in weather forecasting, crop yield prediction, and early warning of flood and drought.
It has been shown that better characterization of soil moisture can significantly improve the weather forecasting. 
In particular, we consider high-resolution daily soil moisture data at the top layer of the 
Mississippi River Basin in the United States, on January 1st, 2004. The spatial resolution is of 
$0.0083$ degrees and the distance of one-degree difference in this region is approximately $87.5$ km.
The grid consists of $1830 \times 1329 = 2{,} 432{,} 070$ 
locations with 2,153,888 measurements and 278,182 missing values. We use the same 
model for the mean process as in Huang and Sun~\cite{huang2017hierarchical}, and fit a zero-mean Gaussian 
process model with a Mat{\'e}rn covariance function to the residuals; see Huang and Sun~\cite{huang2017hierarchical} 
for more details on data description and exploratory data analysis. 

Furthermore, we consider another example of climate and environmental data: wind speed.
Wind speed is an important factor of the climate system's atmospheric quantity. 
It is impacted by the changes in temperature, which lead to air moving from high-pressure
to low pressure layers. Wind speed affects weather forecasting and different
activities related to both air and maritime transportations. Moreover, constructions projects, ranging from
airports, dams, subways and industrial complexes to small housing buildings are impacted by
the wind speed and directions.

The advanced research core of WRF (WRF-ARW) is used in this study to generate
a regional climate dataset over the Arabian Peninsula~\cite{skamarock2005description} 
in the Middle-East. The model is configured
with a domain of a horizontal resolution of $5$ km with $51$ vertical layers while the model top is fixed at $10$ hPa.
The domain covers the longitudes and latitudes of 
$20$\textdegree{}E - $83$\textdegree{}E and $5$\textdegree{}S - $36$\textdegree{}N, respectively. The data are available
daily through $37$ years. Each data file represents $24$ hours measurements of wind speed recorded each hour on
$17$ different layers. In our case, we have picked up one dataset on September 1st, 2017 at time 00:00 AM on a $10$-meter
distance above ground (i.e., layer $0$). No special restriction is applied to the chosen data. 
We only choose an example to show the effectiveness of our proposed framework, but may easily consider 
extending the datasets. 

Since \texttt{ExaGeoStat} can handle large covariance matrix computations, 
and the parallel implementation of the algorithm significantly reduces the 
computational time, we propose to use exact maximum likelihood inference for 
a set of selected regions in the domain of interest to characterize and 
compare the spatial variabilities of both the soil moisture and the wind speed data.


\section{Performance}
\label{sec:perf}
This section evaluates the performance and the accuracy of the 
TLR \texttt{ExaGeoStat} framework for the MLE computations. It presents performance and
accuracy assessments against the reference full accuracy implementation
on shared and distributed-memory systems using synthetic and real datasets.

\subsection{Experimental Settings}
We evaluate the performance of the TLR \texttt{ExaGeoStat} framework for the MLE computations on a wide range of 
Intel hardware systems' generation to highlight our software portability: 
a dual-socket 28-core Intel Skylake Intel Xeon Platinum 8176 CPU running at 2.10 GHz,
a dual-socket 14-core Intel Broadwell Intel Xeon E5-2680 V4 running at 2.4 GHz, 
a dual-socket 18-core Intel Haswell  Intel Xeon CPU E5-2698 v3 running at 2.30 GHz,
Intel manycore Knights Landing (KNL) 7210  chips with 64 cores, 
and a dual-socket 8-core Intel Sandy Bridge Intel Xeon CPU E5-2650 running at 2.00 GHz.
For the distributed-memory experiments, we use \emph{Shaheen-2} from the KAUST Supercomputing Laboratory,
a Cray XC40 system with 6,174 dual-socket compute nodes
based on 16-core Intel Haswell processors running at
2.3 GHz. Each node has 128 GB of DDR4 memory. 
\emph{Shaheen-2} has a total of 197,568
processor cores and 790 TB of aggregate memory. \hatem{In fact, our software portability
is in fact guaranteed, as long as an optimized BLAS/LAPACK high performance library is available
on the targeted system.
}

Our framework is compiled with gcc v5.5.0 and linked against latest
 Chameleon\footnote{https://gitlab.inria.fr/solverstack/chameleon/} 
and HiCMA\footnote{https://github.com/ecrc/hicma} libraries with 
\texttt{HWLOC} v1.11.8, \texttt{StarPU} v1.2.1, \texttt{Intel MKL} v11.3.1, 
\texttt{GSL} v2.4, and \texttt{NLopt} v2.4.2 optimization libraries. All computations
are carried out in double precision arithmetic and each run has been repeated
three times.
The accuracy and qualitative analyses are performed using synthetic and two 
examples of real datasets, i.e., the soil moisture dataset at Mississippi River Basin 
region and the wind speed dataset from the Middle-East region, as described 
in Section~\ref{sec:geospat} for more details.


\subsection{Performance on Shared-Memory Systems}

\begin{figure}
  \centering
  \subfigure[A dual-socket $18$-core Intel Haswell.]{
    \label{fig:time-haswell}
   \includegraphics[width=0.90\linewidth,  height=4.2 cm]{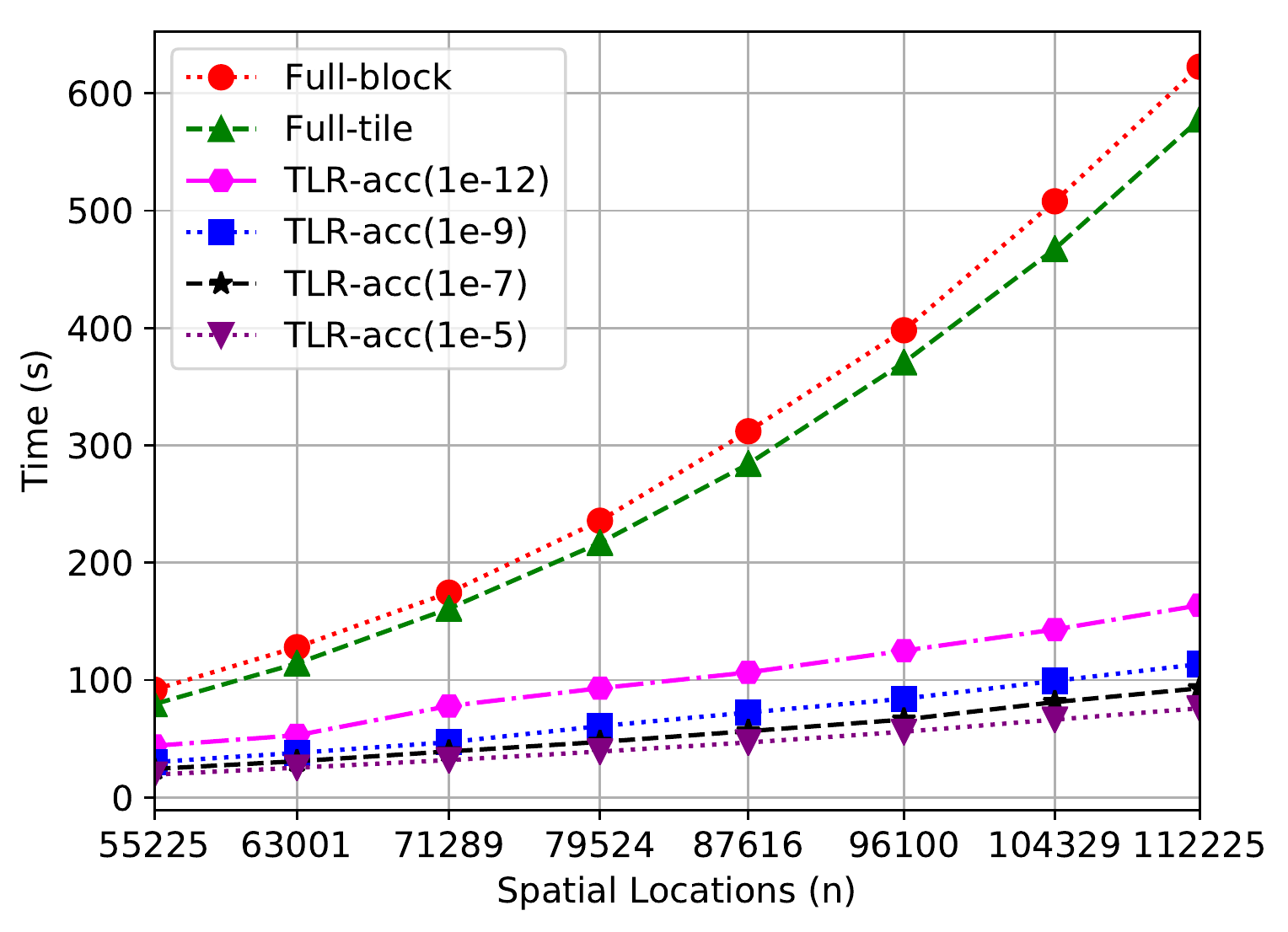}}
  \subfigure[A dual-socket $14$-core Intel Broadwell.]{
    \label{fig:time-broadwell}
   \includegraphics[width=0.90\linewidth,  height=4.2 cm]{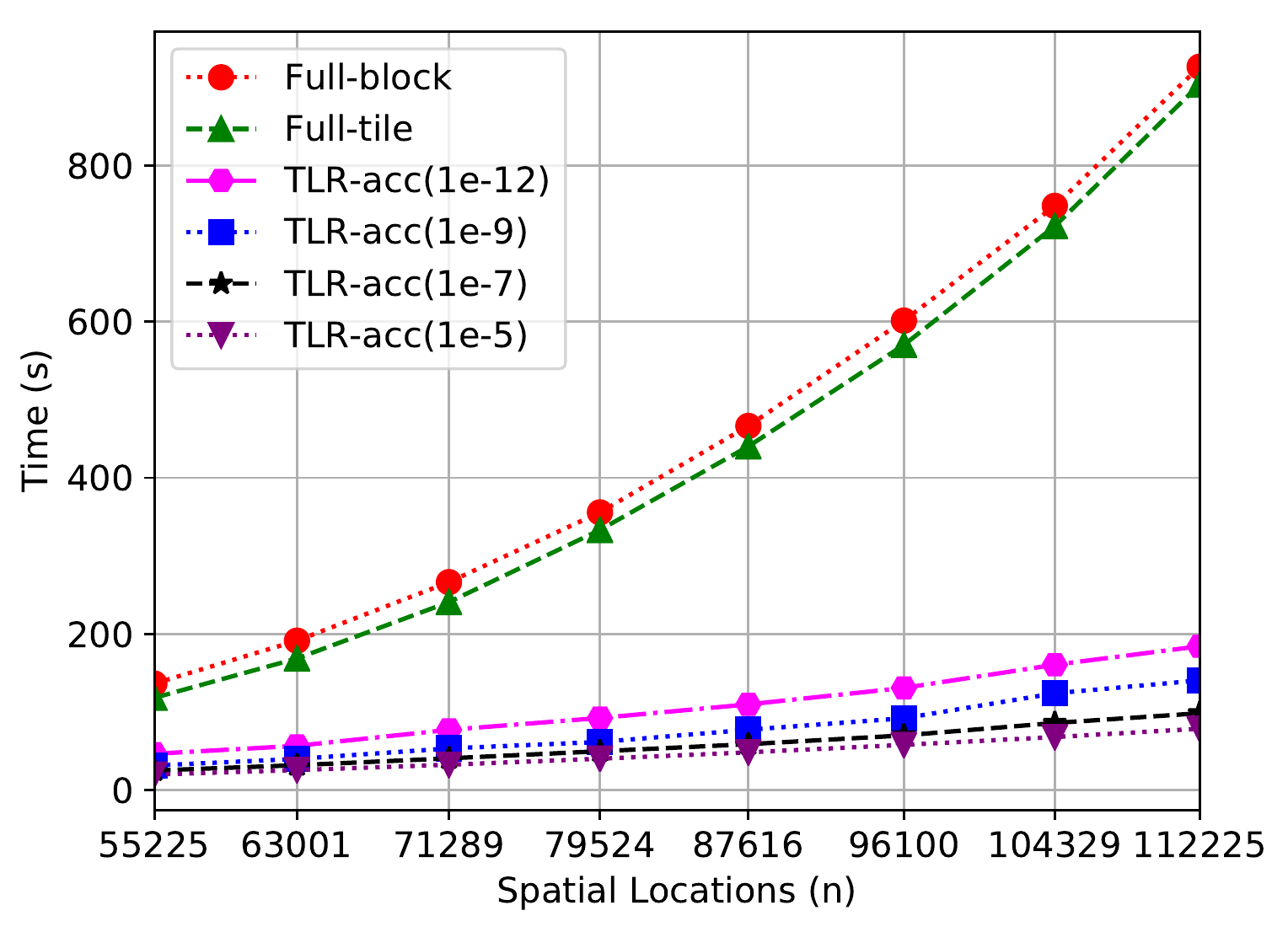}}
  \subfigure[$64$-core Intel Knights Landing (KNL).]{
    \label{fig:time-knl}
   \includegraphics[width=0.90\linewidth,  height=4.2 cm]{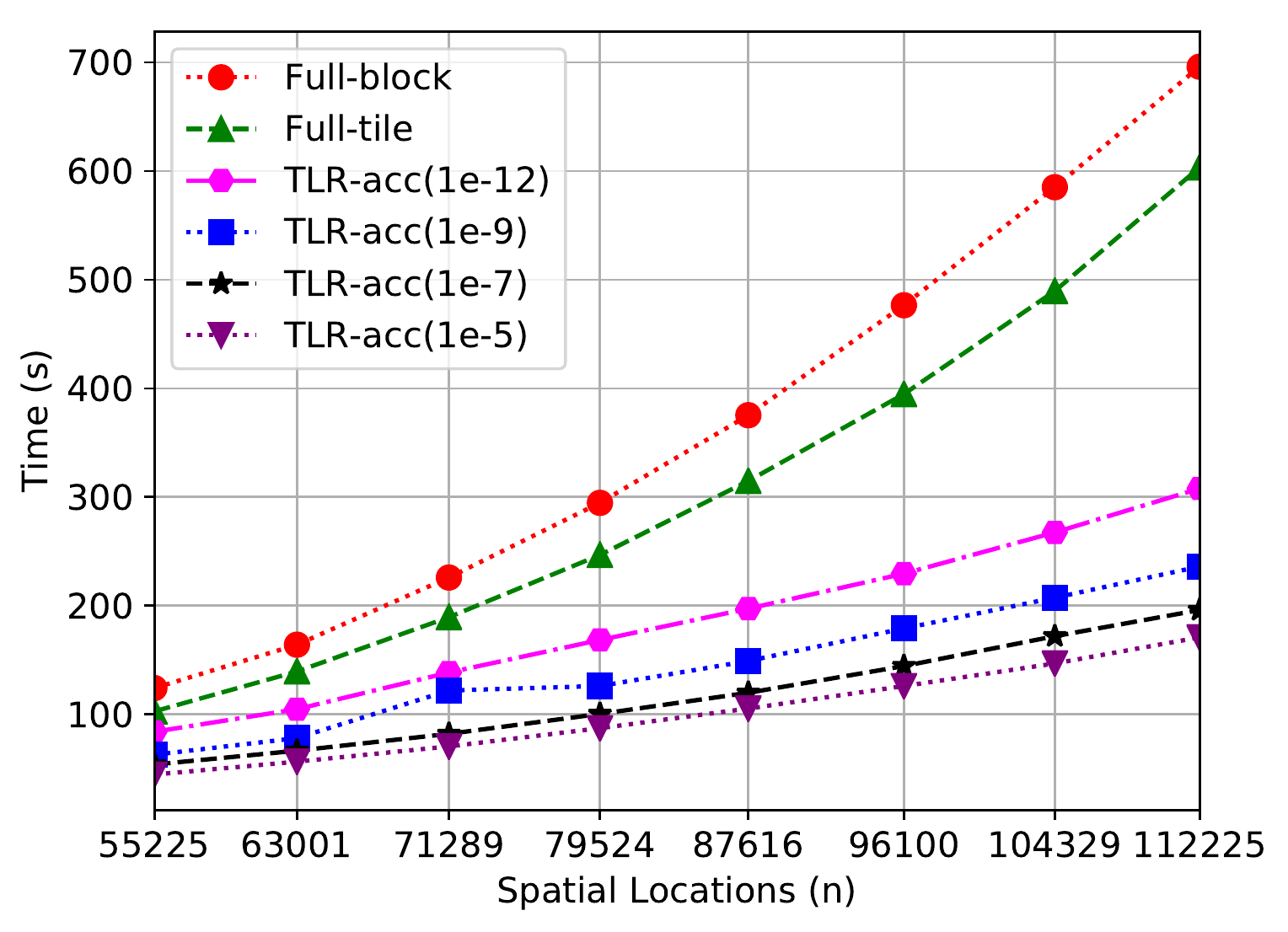}}
 \subfigure[A dual-socket $28$-core Intel Skylake.]{
    \label{fig:time-skylake}
   \includegraphics[width=0.90\linewidth,  height=4.2 cm]{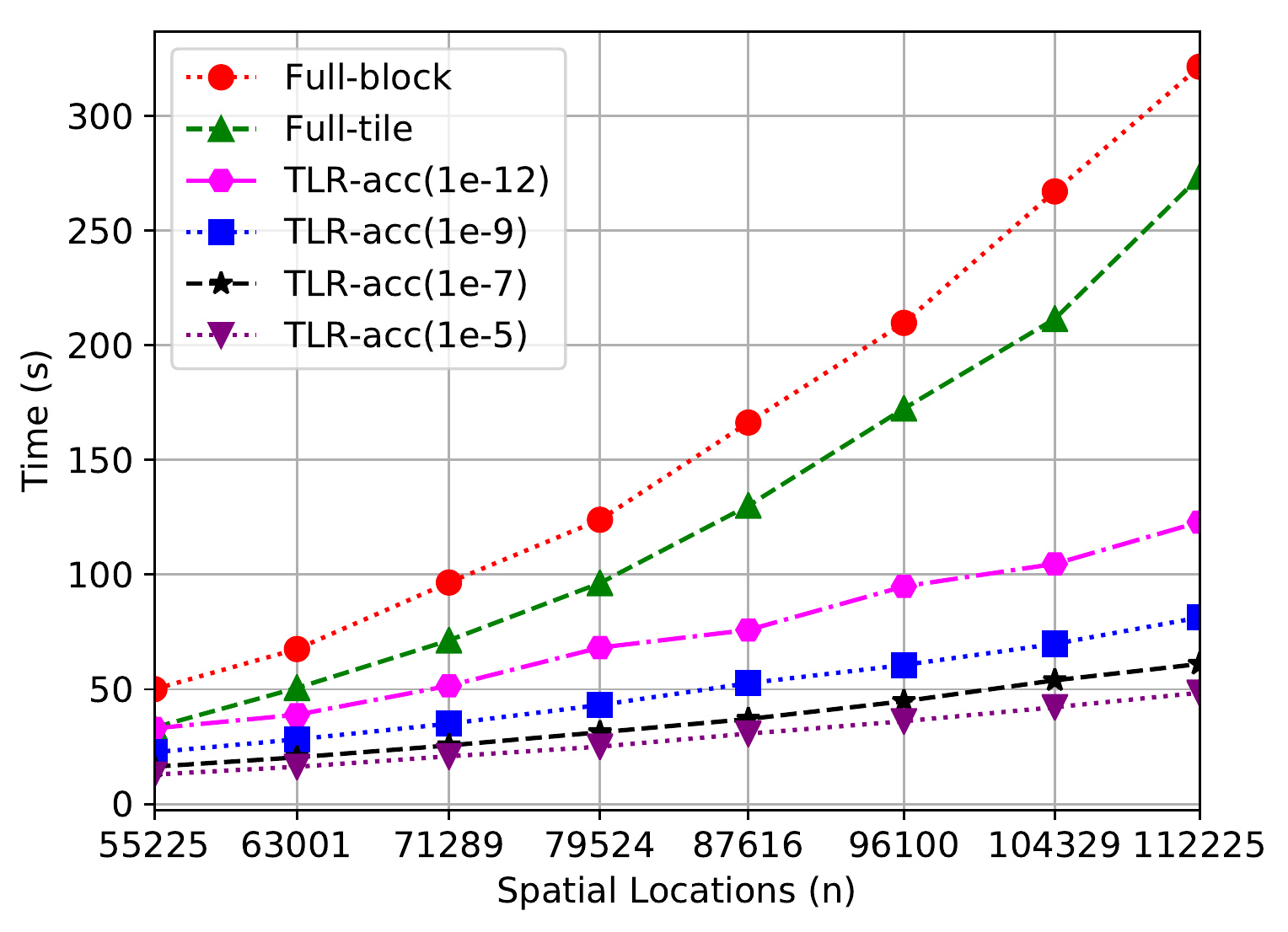}}
   \caption{Time of one iteration of the TLR MLE operation on different Intel architectures.}
   \label{fig:perf_comp}
\end{figure}

We present the performance analysis 
of TLR MLE computation on the four aforementioned Intel systems over various
numbers of spatial locations. We compare  
against the full machine precision accuracy obtained from the block and tile MLE implementations,
with the \texttt{Intel MKL LAPACK} and \texttt{Chameleon} libraries, respectively, as described in 
Section~\ref{sec:tlr}. 
In the following figures, the x-axis represents the number of spatial locations, and
the y-axis represents the total execution time in seconds. We use 
four different TLR accuracy thresholds, i.e., $10^{-5}$, $10^{-7}$, $10^{-9}$, and $10^{-12}$.
%
Figure~\ref{fig:perf_comp} shows the time to solution to perform 
the TLR MLE operation across generations of Intel systems. Since the internal optimization process
is an iterative procedure, which usually takes a few tens of iterations, 
we only report the time for a single iteration as a proxy for the overall optimization procedure. 
The elapsed time for the full machine precision accuracy for the tile MLE outperforms the block
implementations. This is expected and has already been highlighted in~\cite{abdulah2017exageostat}.
Regarding the TLR MLE implementations, the time to solution
steadily diminishes as the requested accuracy decreases. This phenomenon is reproduced
on all shared-memory platforms. The maximum speedup achieved by TLR MLE is significant for 
all studied accuracy thresholds. In particular, the maximum speedup obtained with $10^{-5}$
accuracy threshold is around $7$X, $10$X, $13$X and $5$X on the Intel Haswell, 
Broadwell, KNL and Skylake, respectively.
The  speedup numbers have to be cautiously assessed, since approximation
obviously introduces numerical errors and may not be bearable by the geospatial statistics application
beyond a certain threshold. So, the challenge is to maintain the model fidelity with
high accuracy, just enough, to outperform the full machine precision accuracy 
MLE implementations by a non-negligible factor.

 

   \begin{figure}
  \centering
  \subfigure[$256$ nodes.]{
    \label{fig:time-shaheen-256}
   \includegraphics[width=0.90\linewidth,  height=4.2 cm]{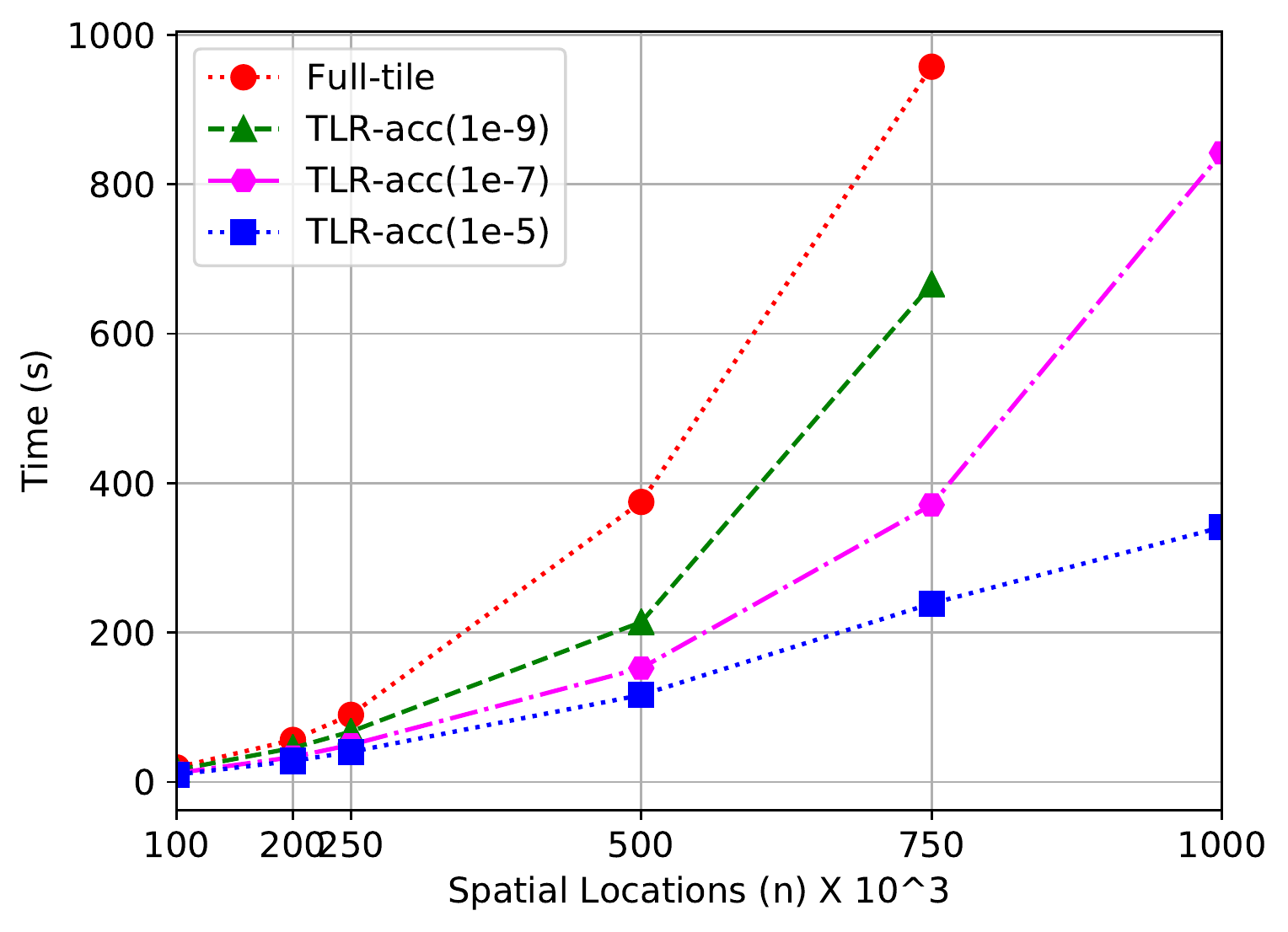}}
     \subfigure[$1024$ nodes.]{
    \label{fig:time-shaheen-1024}
   \includegraphics[width=0.90\linewidth,  height=4.2 cm]{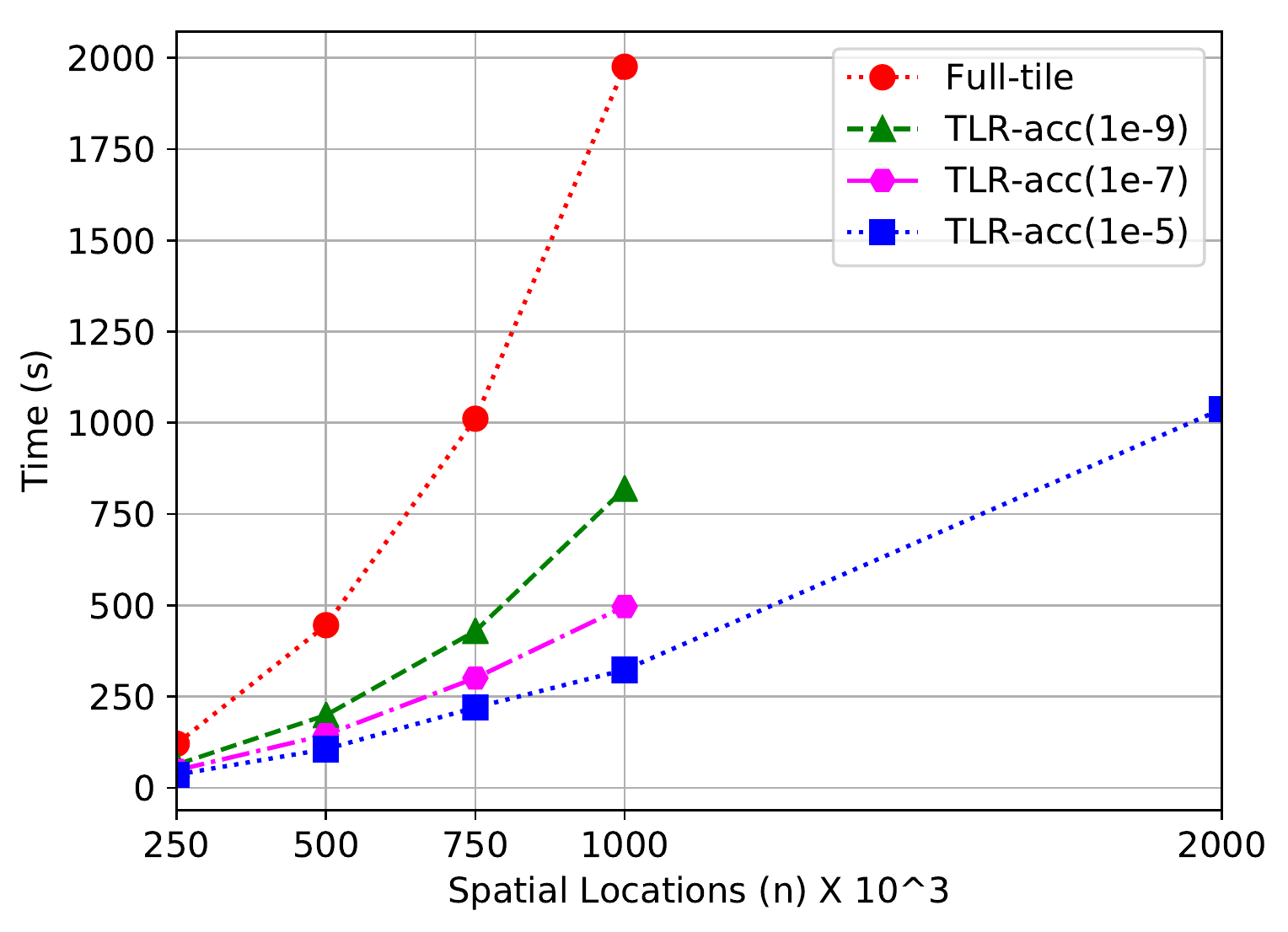}}
   \caption{Time of one iteration of the TLR MLE operation on Cray XC40 \emph{Shaheen-2} using different accuracy thresholds.}
\label{fig:time-shaheen}
\end{figure}


%
%
%
%


\subsection{Performance on Distributed-Memory Systems}

We also test our TLR \texttt{ExaGeoStat} framework on the distributed-memory 
\emph{Shaheen-2} Cray XC40 system using $256$ ($\sim 8,200$ cores) and $1024$ nodes ($\sim 33,000$ cores), as highlighted in 
Figure~\ref{fig:time-shaheen}. Similarly to the results on shared-memory systems,
significant speedups (up to $5$X) are achieved when performing TLR \texttt{ExaGeoStat} approximations for the MLE
computations using different numbers of spatial locations. There are some points
missing in the graphs, which correspond to cases where the application runs out of memory.
While this is the case when performing no TLR approximation, this may also occur with TLR 
approximation with high accuracy thresholds.
Tuning the tile size $nb$ is of paramount importance to achieving high performance when running in approximation 
or full accuracy mode. For instance, for the full machine precision accuracy (i.e, \emph{full-tile}) 
variant of the \texttt{ExaGeoStat} MLE calculations,  
we use a tile size of $560$, while a much higher tile size of $1900$ is required for the TLR variants to
be competitive. This tile size discrepancy is due to the resulting arithmetic intensity of the main
computational tasks for each variant. For the \emph{full-tile} variant, since the main kernel is
the dense matrix-matrix multiplication, $nb=560$ is large enough to keep single cores busy caching data
located at the high level of the memory subsystem, and small enough to provide enough concurrency
for the overall execution of the MLE application.
For the TLR variants, large tile size is necessary, since the shape of the data structure depends
on the actual rank obtained after compression. These ranks are usually much smaller than the tile 
size $nb$. Moreover, the main computational kernel for TLR MLE computations is the TLR 
matrix-matrix multiplication, which involves several successive linear algebra calls~\cite{akbudak-isc17}.
As a result, the resulting arithmetic intensity of that kernel is rather low, close to memory-bound regime.
In distributed-memory systems, this mode translates into latency-bound since data motion happens
between remote node memories. This engenders significant overheads, which can not be compensated since 
computation is very limited.
We may therefore increase the tile size just enough to slightly shift the 
regime toward compute-bound, while preserving a high level of parallelism to maintain hardware occupancy. 
We tuned the tile size $nb$ on our target distributed-memory  \emph{Shaheen-2} Cray XC40 system to gain the best performance of our MLE implementation.

\begin{figure}
\centering
\includegraphics[width=0.90\linewidth,  height=4 cm]{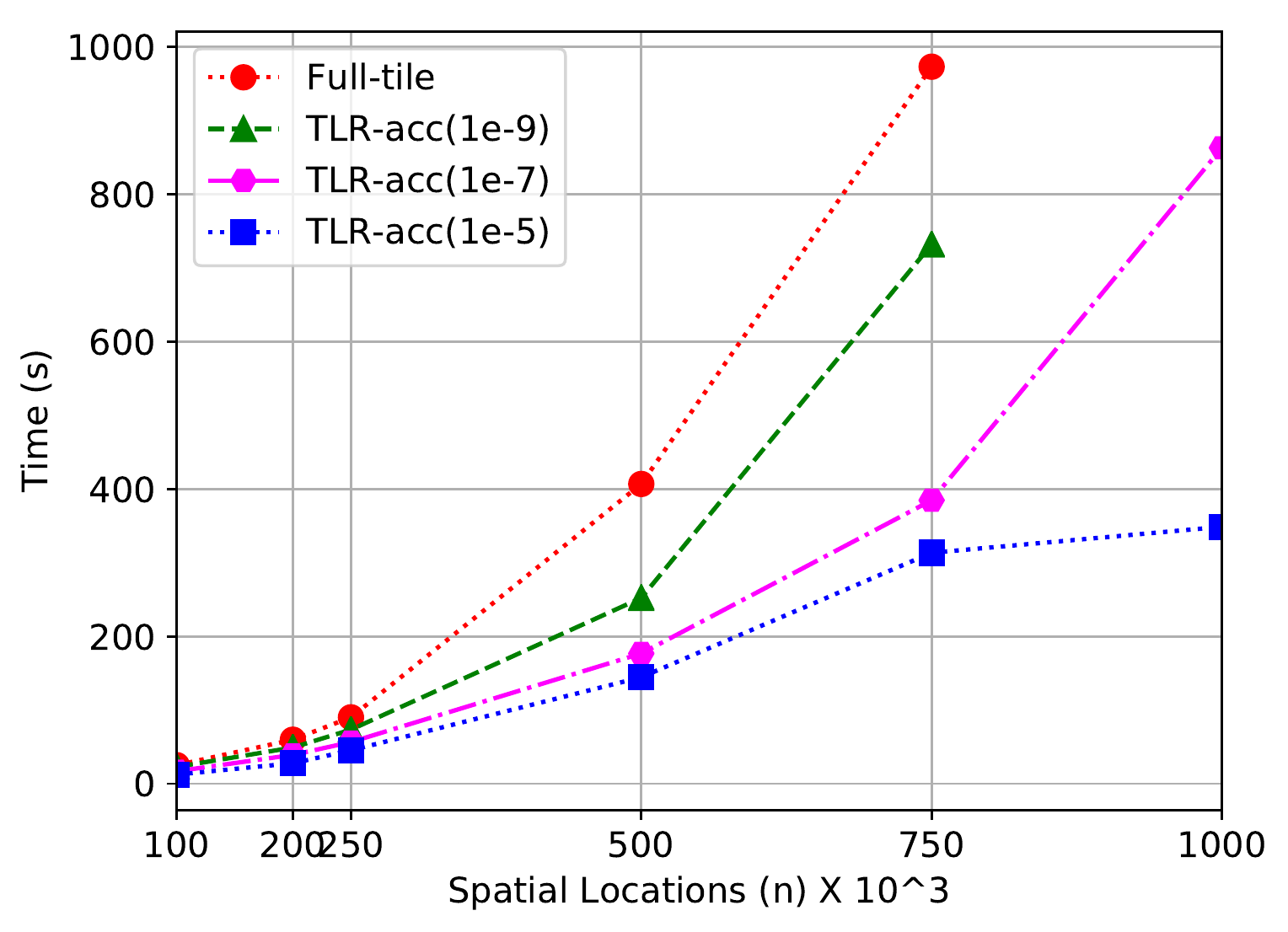}
  \caption{Time of TLR prediction operation on Cray XC40 \emph{Shaheen-2} with $256$ nodes.}
  \label{fig:time-pred-shaheen}
\end{figure}

Moreover, we  investigate the performance of the prediction operation (i.e., $100$ unknown measurements),
as introduced in equation (\ref{eq:pred3}). 
Figure~\ref{fig:time-pred-shaheen} shows the execution time for the prediction
  on different synthetic datasets up to $1 M \times 1 M$ matrix size using $256$ nodes.
   The most time-consuming part of the prediction operation is the Cholesky factorization in this configuration, 
   since the number of unknown measurements to calculate is rather small and triggers only a small number of triangular
   solves.
    Thus, the performance curves show a similar behavior as the MLE operation using the same
     number of nodes, as shown in Figure~\ref{fig:time-shaheen-256}.

\subsection{Accuracy Verification and Qualitative Analysis}
\label{sec:acc-comp}
We evaluate the accuracy
of TLR approximation techniques for the MLE calculations with different accuracy
 thresholds and compare against its full machine precision accuracy variant.
The accuracy can be verified at two different occasions: estimating the 
MLE parameter vector and predicting missing
measurements at certain locations. Here, we 
use both synthetic and real datasets to perform this
accuracy verification and analyze the effectiveness
of the proposed approximation techniques. 

\subsubsection{Synthetic Datasets (Monte Carlo Simulations)}

\begin{figure*}[!ht]
  \centering
  \subfigure[Estimated variance parameter ($\theta_1$).]{
    \label{fig:theta1-1-0.03-0.5}
   \includegraphics[width=0.31\linewidth]{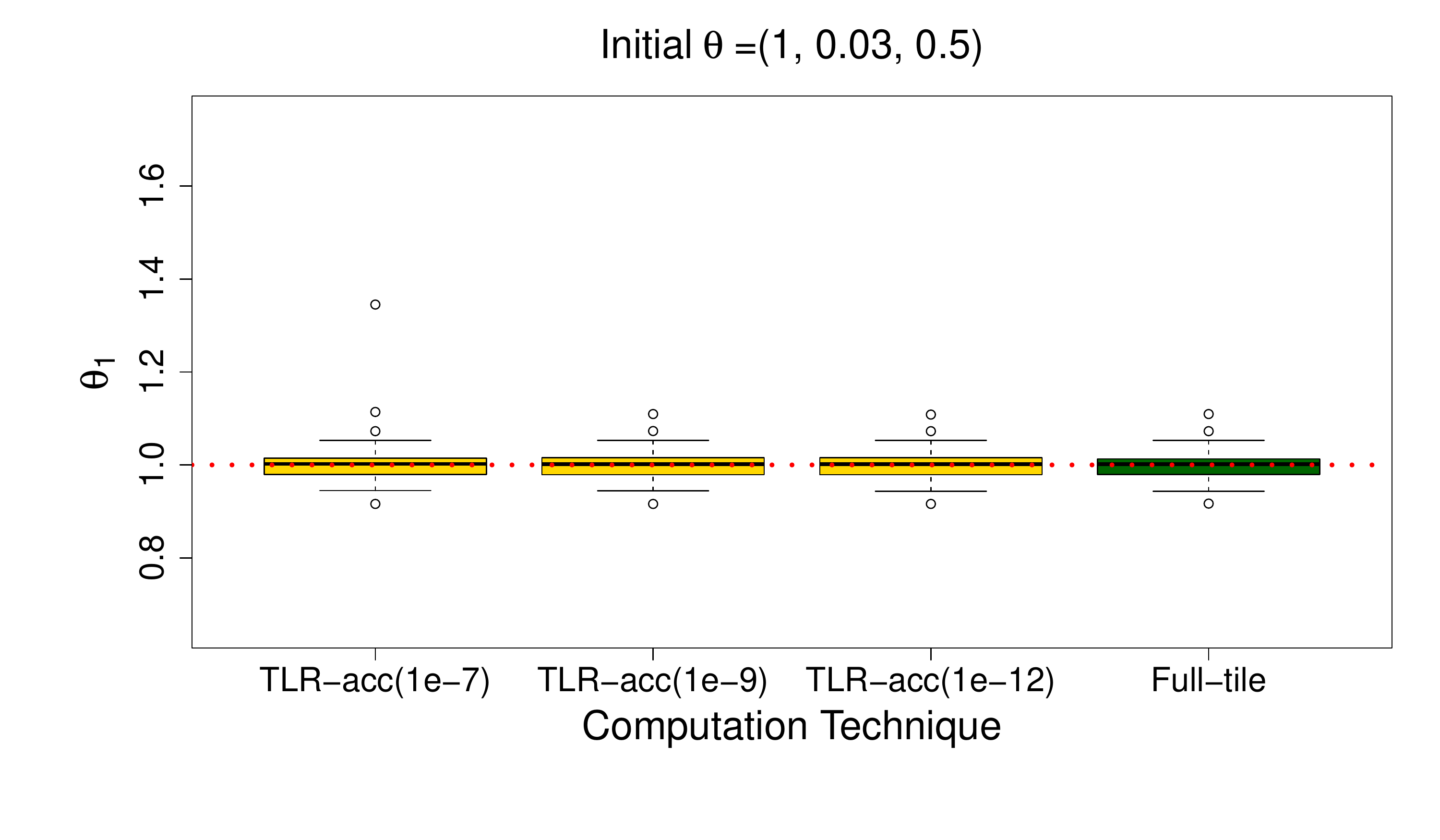}}
  \subfigure[Estimated spatial range parameter ($\theta_2$).]{
    \label{fig:theta2-1-0.03-0.5}
   \includegraphics[width=0.31\linewidth]{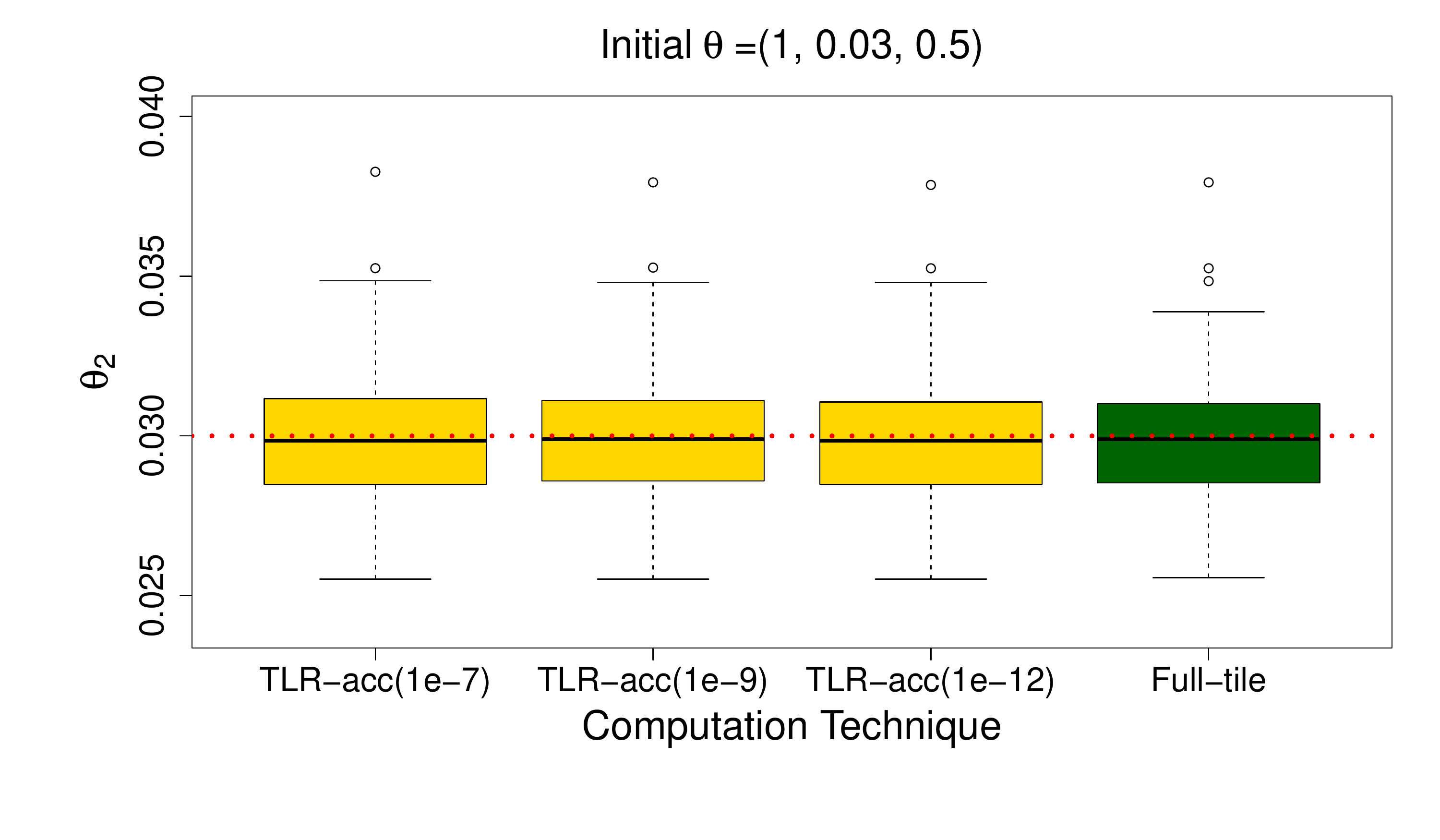}}
  \subfigure[Estimated smoothness parameter ($\theta_3$).]{
    \label{fig:theta3-1-0.03-0.5}
   \includegraphics[width=0.31\linewidth]{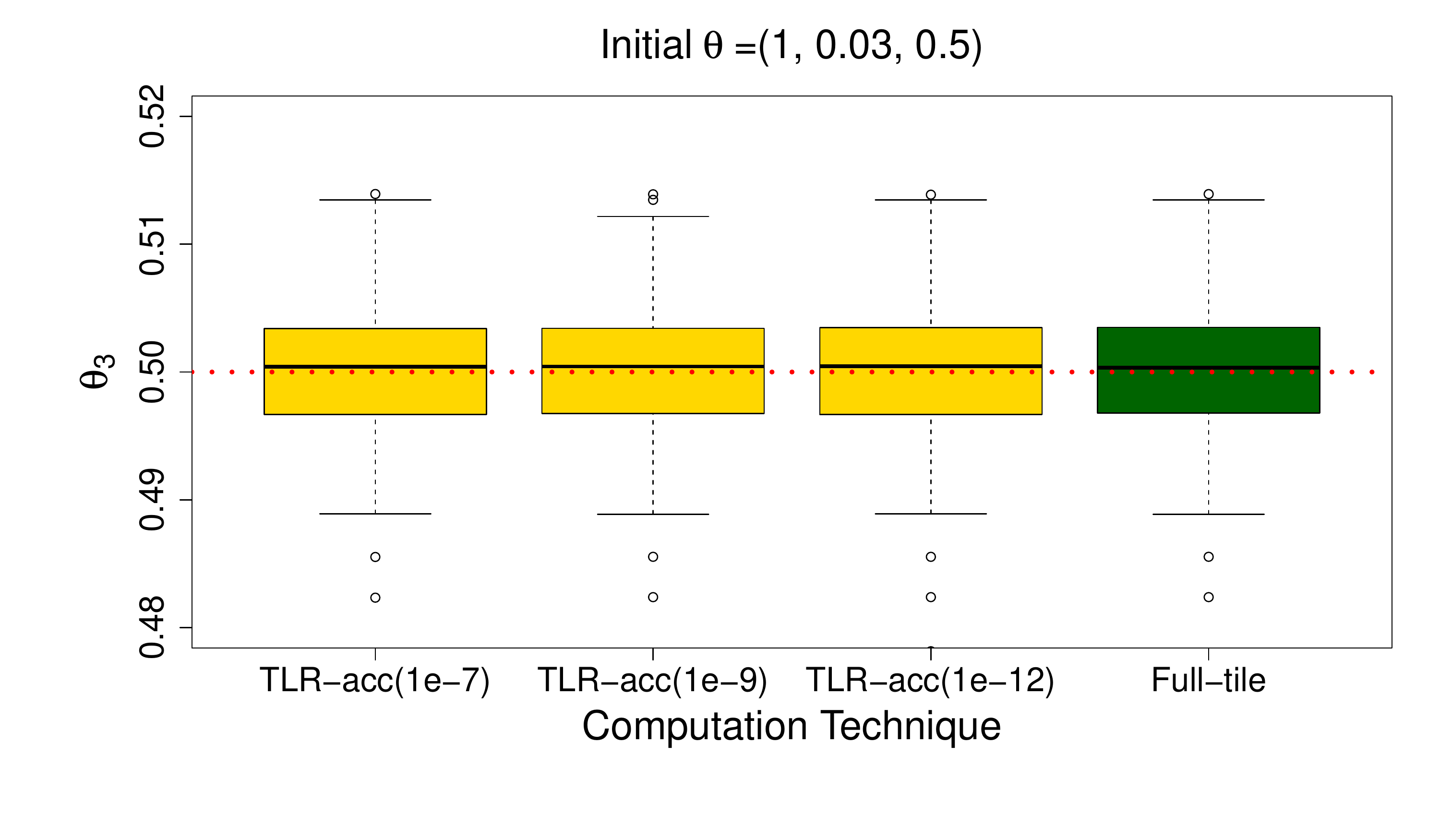}}
  \subfigure[Estimated variance parameter ($\theta_1$).]{
    \label{fig:theta1-1-0.1-0.5}
   \includegraphics[width=0.31\linewidth]{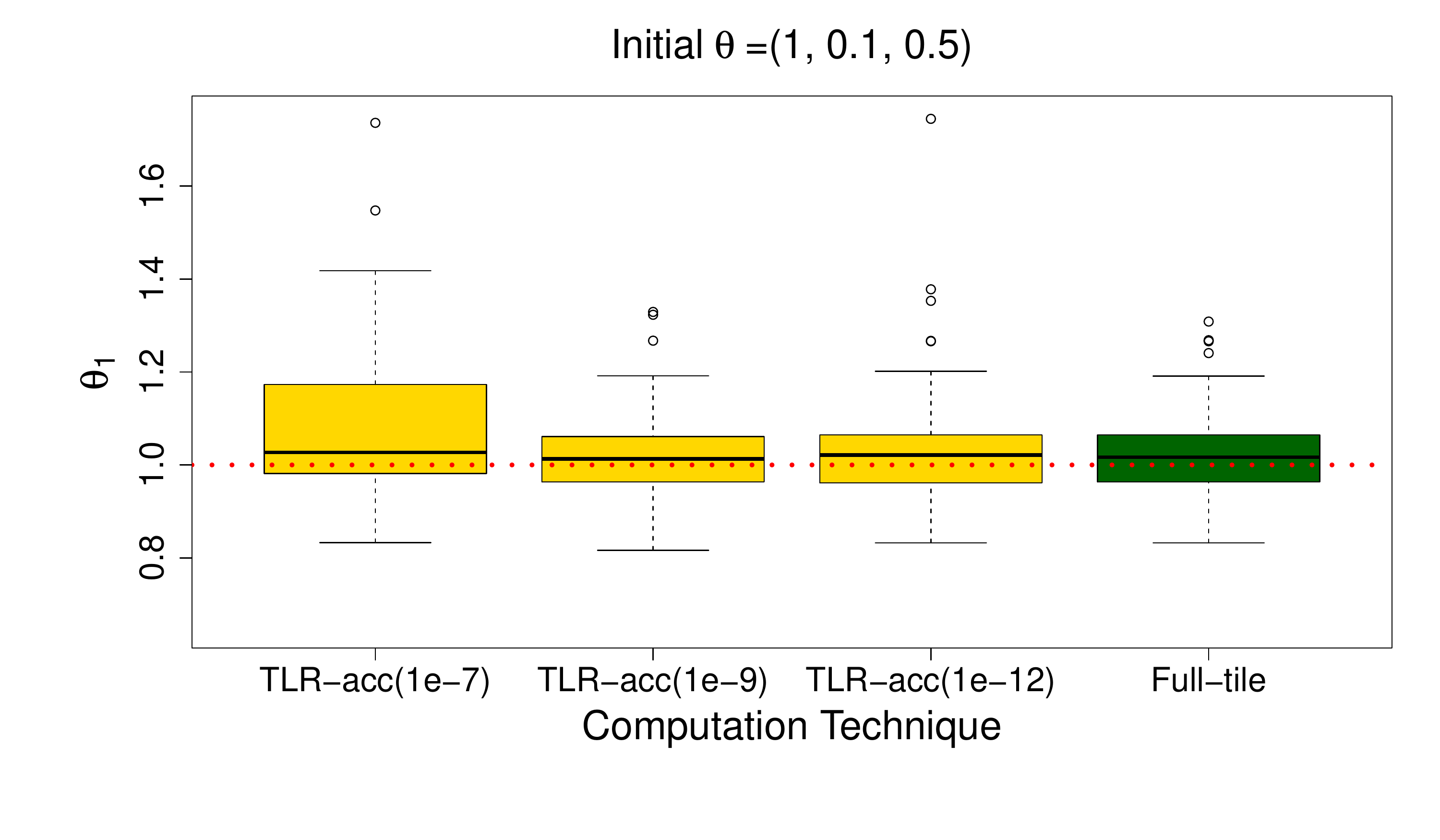}}
  \subfigure[Estimated spatial range parameter ($\theta_2$).]{
    \label{fig:theta2-1-0.1-0.5}
   \includegraphics[width=0.31\linewidth]{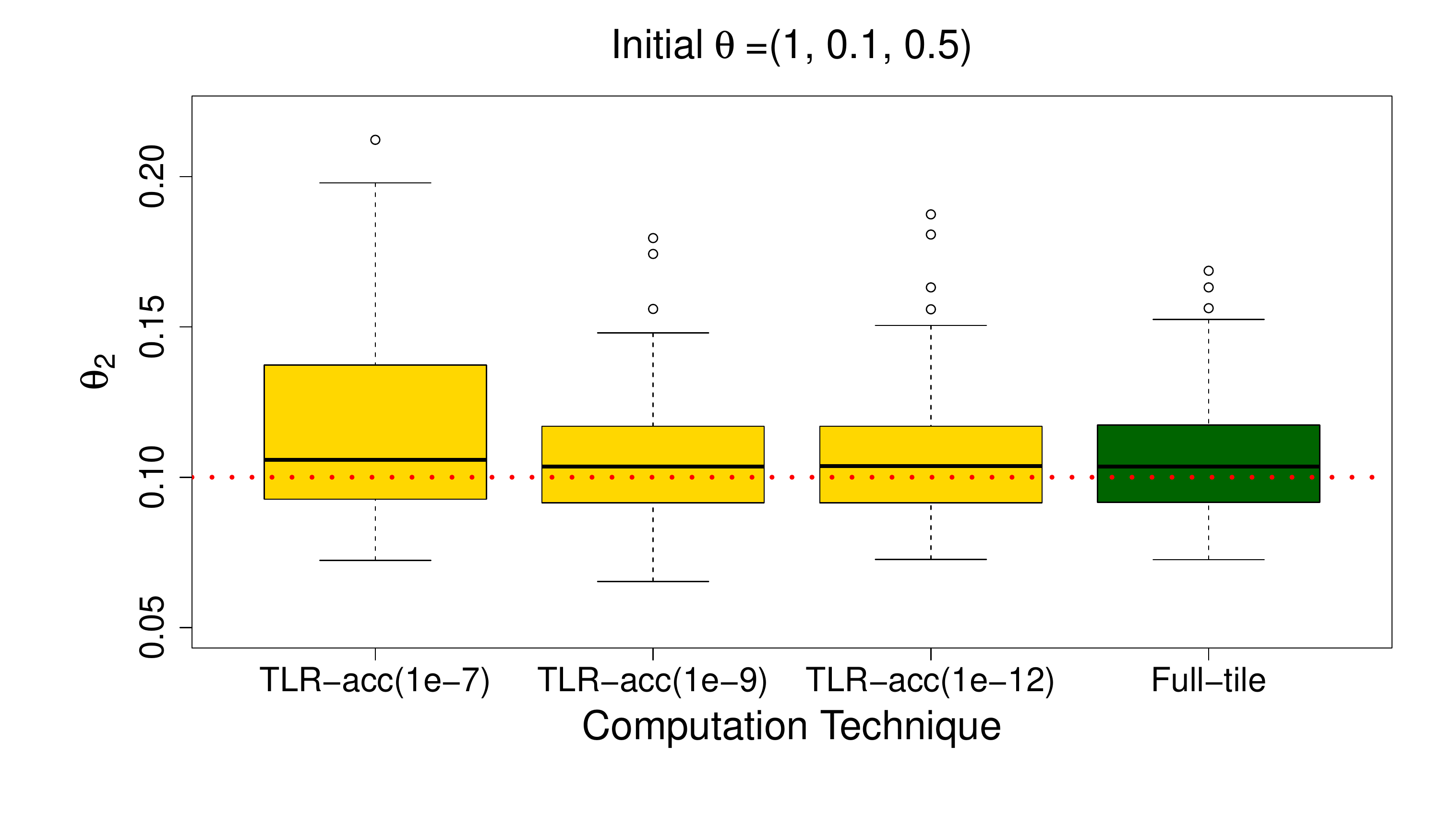}}
  \subfigure[Estimated smoothness parameter ($\theta_3$).]{
    \label{fig:theta3-1-0.1-0.5}
   \includegraphics[width=0.31\linewidth]{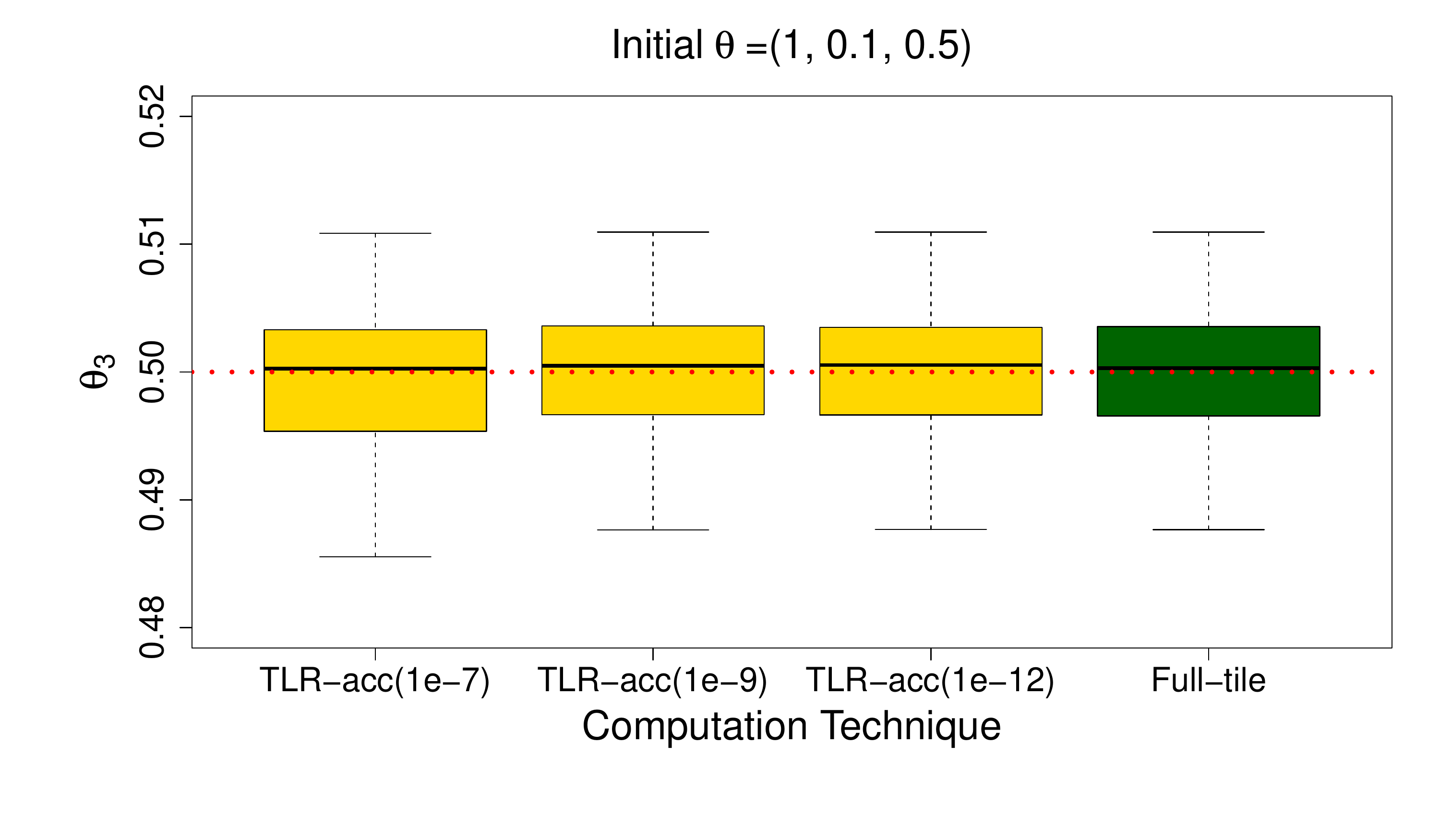}}
     \subfigure[Estimated variance parameter ($\theta_1$).]{
    \label{fig:theta1-1-0.3-0.5}
   \includegraphics[width=0.31\linewidth]{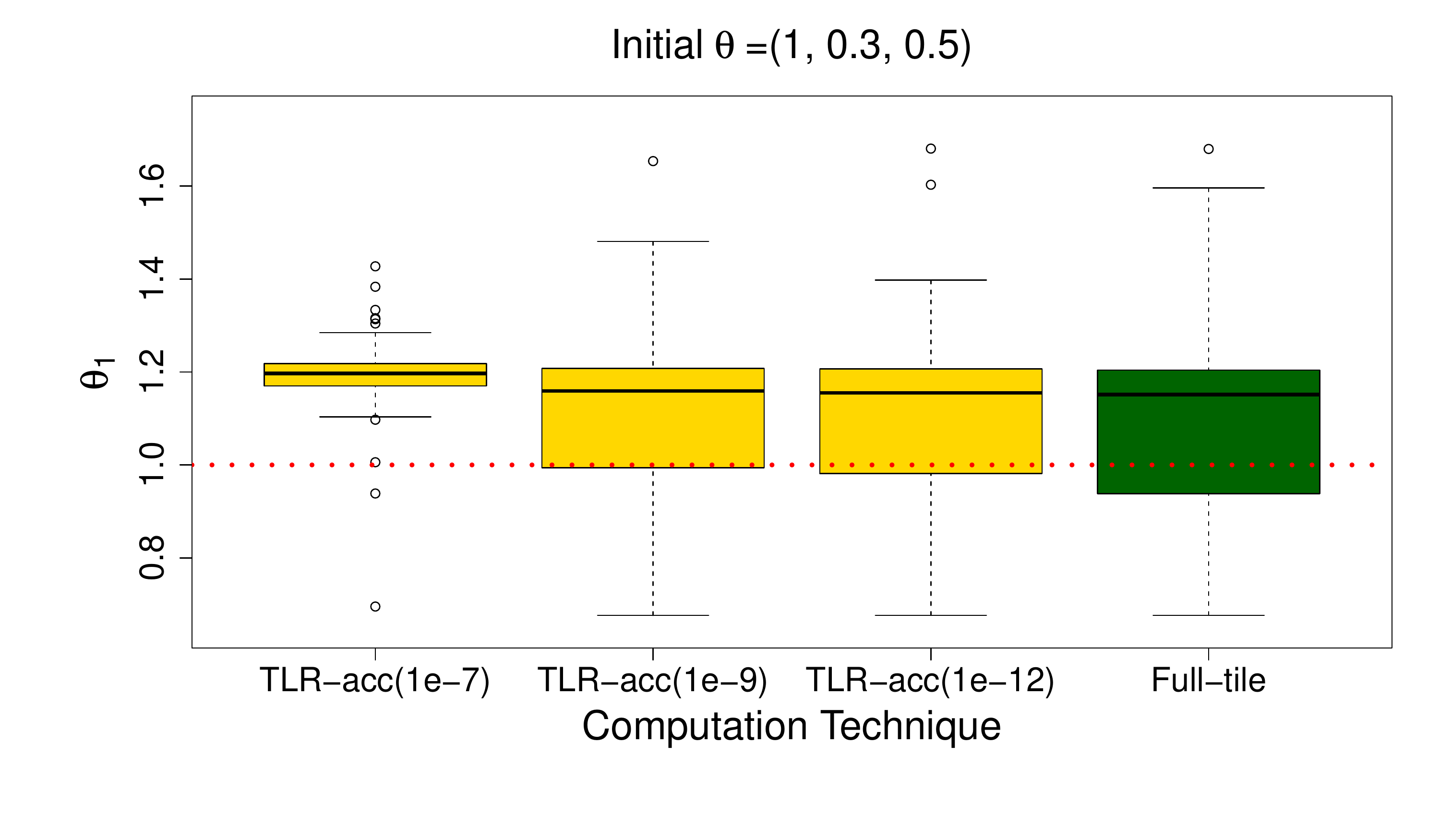}}
  \subfigure[Estimated spatial range parameter ($\theta_2$).]{
    \label{fig:theta2-1-0.3-0.5}
   \includegraphics[width=0.31\linewidth]{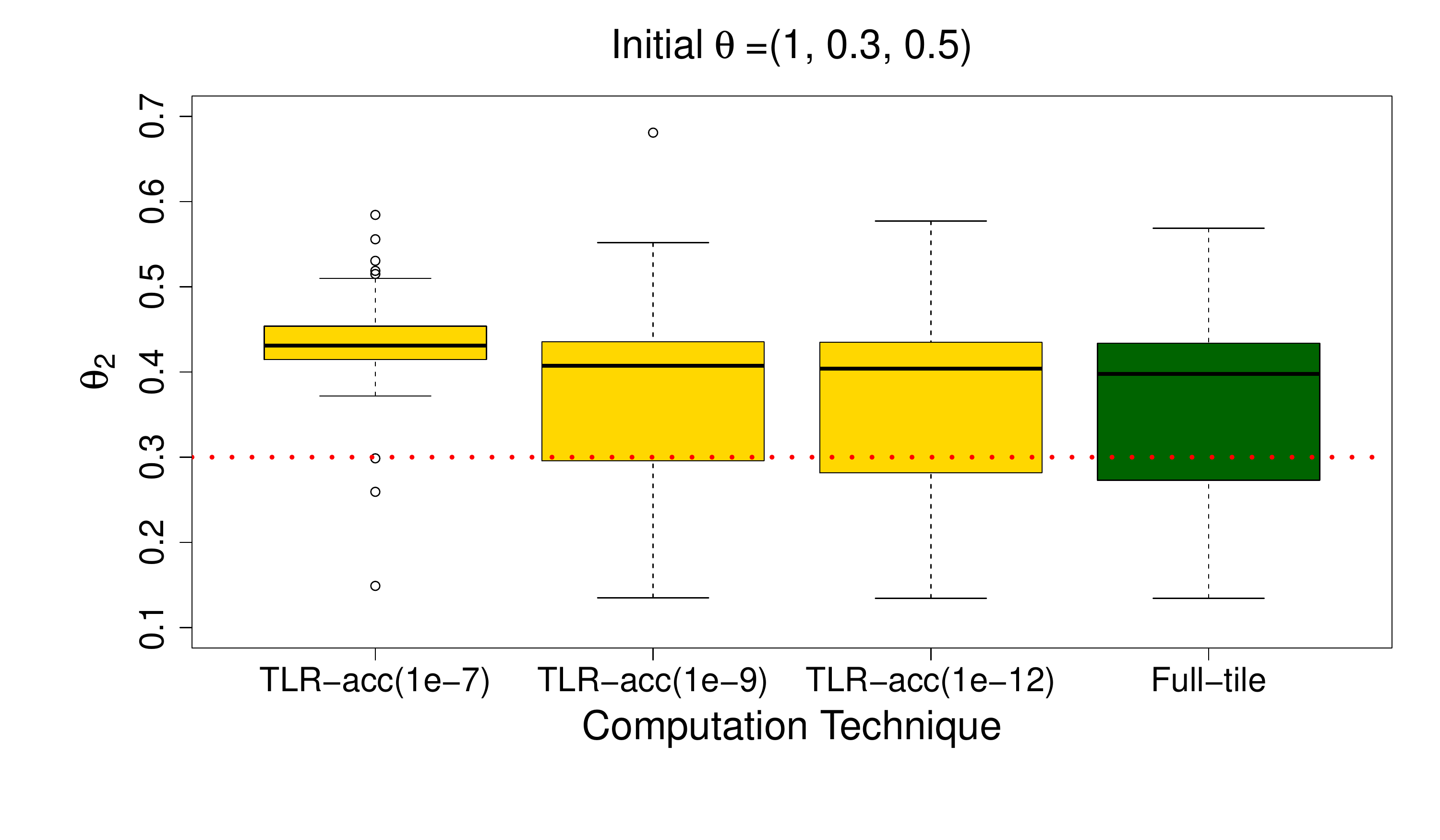}}
  \subfigure[Estimated smoothness parameter ($\theta_3$).]{
    \label{fig:theta3-1-0.3-0.5}
   \includegraphics[width=0.31\linewidth]{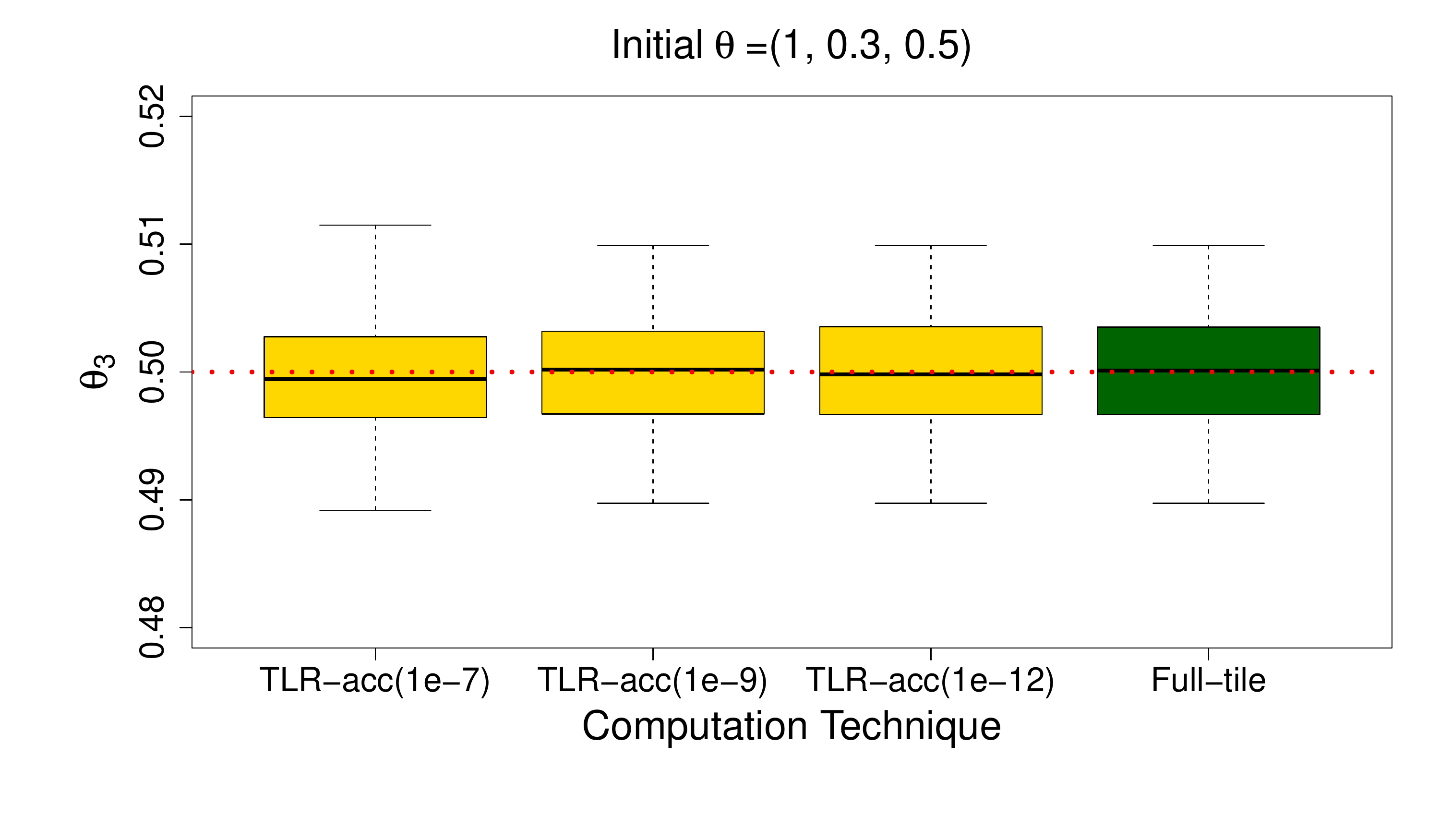}}
\caption{Boxplots of parameter estimation ($\theta_1$, $\theta_2$, and $\theta_3$). 
}
\label{fig:dep_comp}
\end{figure*}

The overall goal of the MLE model is to estimate the unknown
parameters of the underlying statistical model $(\theta_1, \theta_2, \theta_3)$ 
of the Mat{\'e}rn covariance function, then to use this model
for future predictions.  Monte Carlo simulation is a common
way to estimate
the accuracy of the MLE operation using synthetic datasets. 
{\texttt ExaGeoStat} data generator is used to generate
 synthetic datasets. The input of the data generator is an
initial parameter vector to produce a set of locations and
measurements. This initial parameter vector can be reproduced
by the MLE operation using the generated spatial data. 
More details about the data generation process can be found
 in~\cite{abdulah2017exageostat}. We  generate a $40 K$ 
 synthetic data, one location matrix and 100 different 
 measurement vectors,  in exact computation. We rely on exact 
 computation on this step to ensure that all techniques are using the same data for the MLE operation.

As described in Section~\ref{sec:matern_kernel}, the Mat\'{e}rn covariance
function depends on three parameters, variance $\theta_1$, range $\theta_2$, and smoothness $\theta_3$.
The correlation strength can be determined using the range parameter 
$\theta_2$ (i.e., strong correlated data ($\theta_2=0.3$),  
medium correlated data ($\theta_2=0.1$), weak correlated data ($\theta_2=0.03$). 
These correlations values are restricted by the smoothness parameter ($\theta_3=0.5$).
Thus, we select three combination of the three parameters, i.e., 
($(1, 0.3, 0.5)$, $(1, 0.1, 0.5)$, and $(1, 0.03, 0.5)$). The correlation
has obviously a direct impact on the compression rate, and therefore, the actual ranks
of the TLR covariance matrix for the MLE computation.

Figure~\ref{fig:dep_comp} shows three boxplots for each initial parameter vector representing
 the estimation accuracy using different computation techniques.  The true value of each
  $\theta$ is denoted by a dotted red line. It is clearly seen from the figure that
   with weak correlation (i.e., $\theta_2 = 0.03$), TLR approximation is able to retrieve
    the initial parameter vector with the same accuracy as the exact computation. 

The medium correlation (i.e., $\theta_2 = 0.1$) also shows better accuracy with TLR
 approximation up to accuracy $10^{-9}$. However, TLR with accuracy $10^{-7}$
  is less compared to other accuracy levels. With stronger correlation (i.e., $\theta_2 = 0.3$),
   TLR with accuracy $10^{-7}$ and $10^{-9}$
 are not able to retrieve the parameter vector efficiently. In this case, only TLR with
  accuracy $10^{-12}$ can be compared with the exact solution. In summary,
   TLR requires higher accuracy if the data is strongly correlated.



Prediction is key to checking the TLR approximation accuracy
 compared to the \emph{full-tile} variant.
Here, we conduct another experiment to predict $100$ missing values from synthetic
 datasets generated from our three parameter vectors (i.e., $(1, 0.03, 0.5)$, $(1, 0.1, 0.5)$, and $(1, 0.3, 0.5)$).
  The missing values are randomly picked from the generated data so that it can be
   used as a prediction accuracy reference. To assess the accuracy, we use the Mean Square Error (MSE) metric as follows:
\begin{equation}
	\label{eq:MSE}
{\displaystyle \operatorname {MSE} ={\frac {1}{100}}\sum _{i=1}^{100}(Y_{i}-{\hat {Y_{i}}})^{2}.}
\end{equation}
\begin{figure*}[!ht]
  \centering
  \subfigure[Initial $\boldsymbol \theta$ vector $(1-0.03-0.5)$.]{
    \label{fig:theta1p-1-0.1-0.5}
   \includegraphics[width=0.31\linewidth]{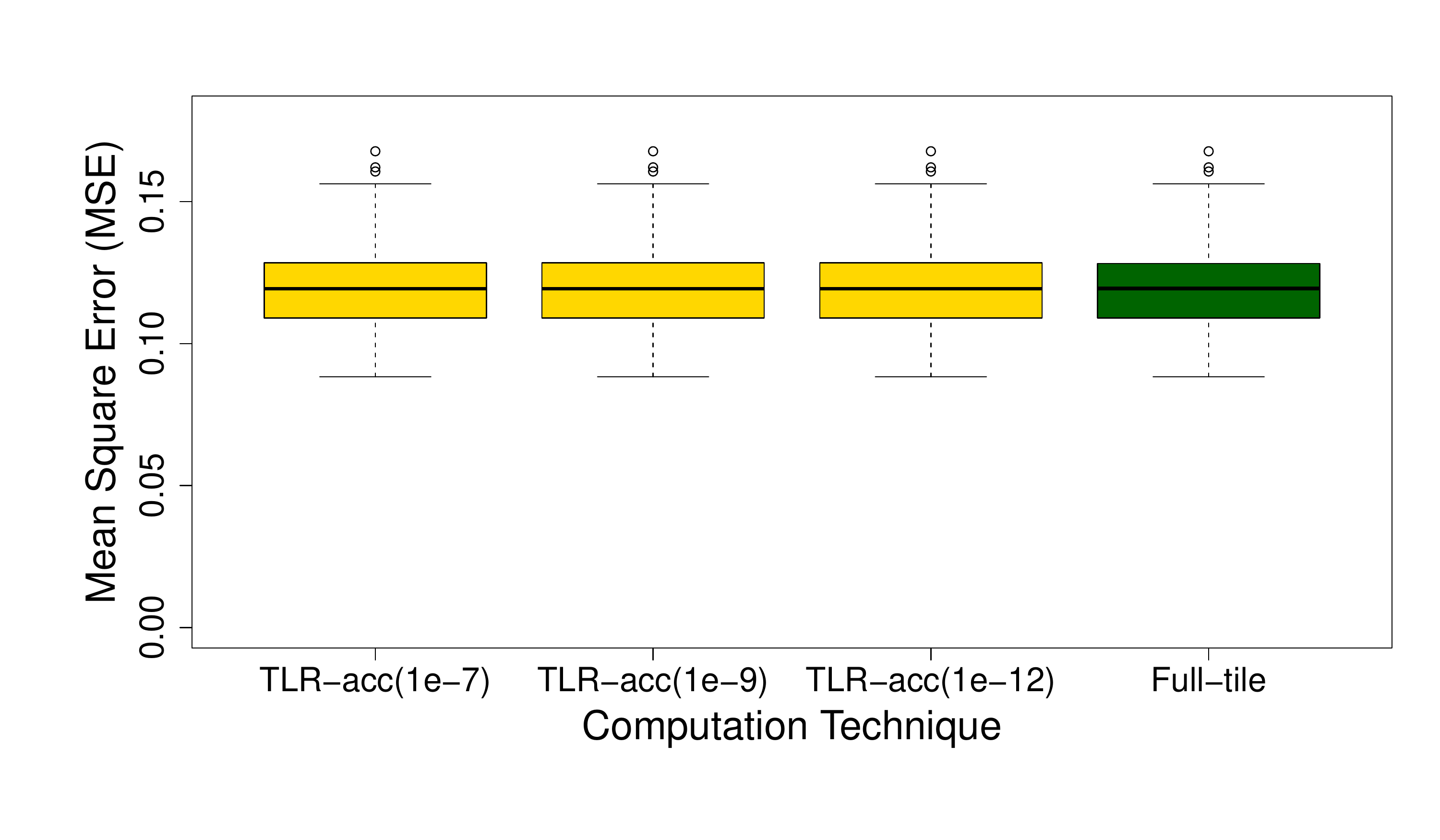}}
  \subfigure[Initial $\boldsymbol \theta$ vector $(1-0.1-0.5)$.]{
    \label{fig:theta2p-1-0.03-0.5}
   \includegraphics[width=0.31\linewidth]{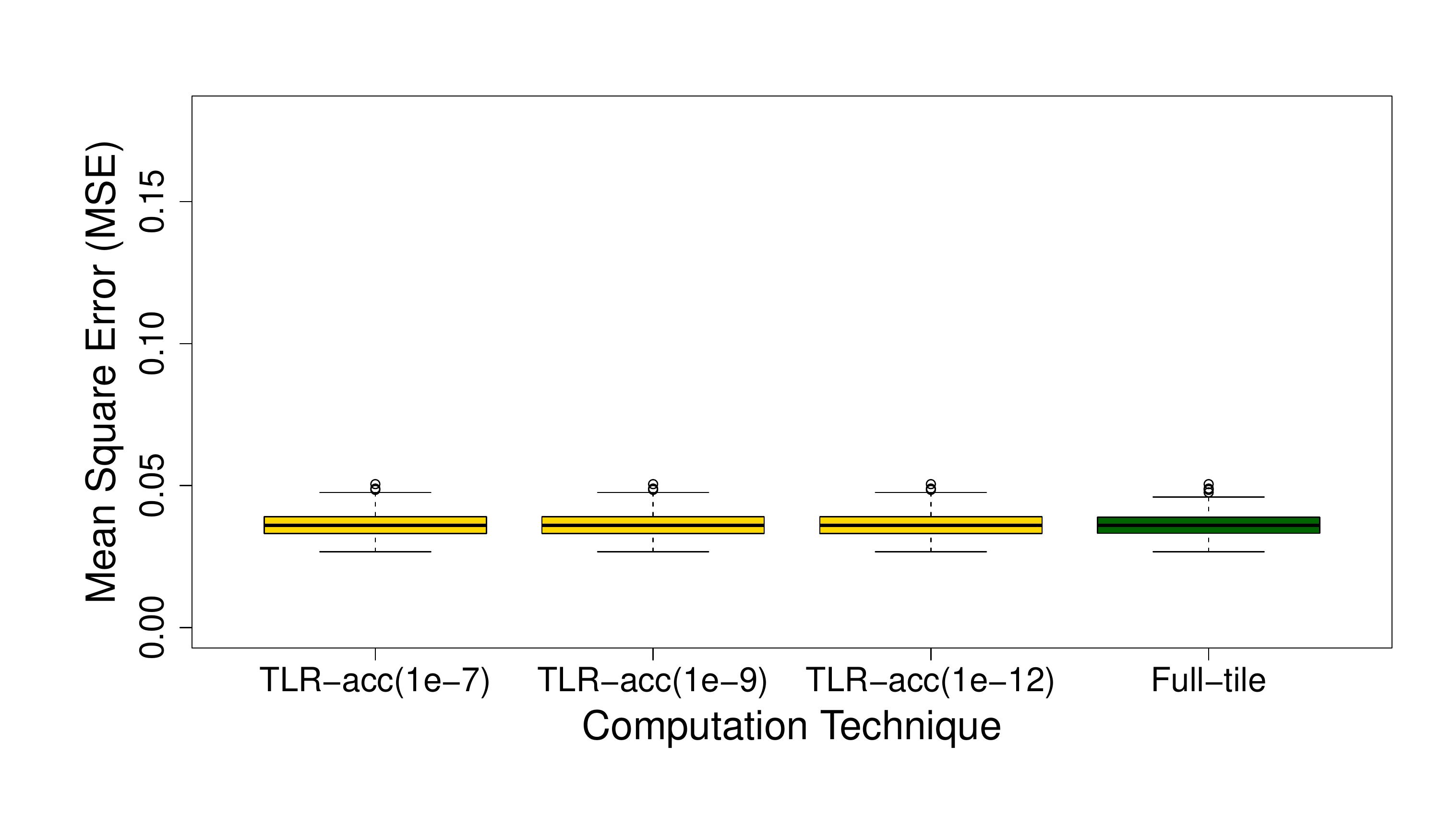}}
  \subfigure[Initial $\boldsymbol \theta$ vector $(1-0.3-0.5$).]{
    \label{fig:theta3p-1-0.3-0.5}
   \includegraphics[width=0.31\linewidth]{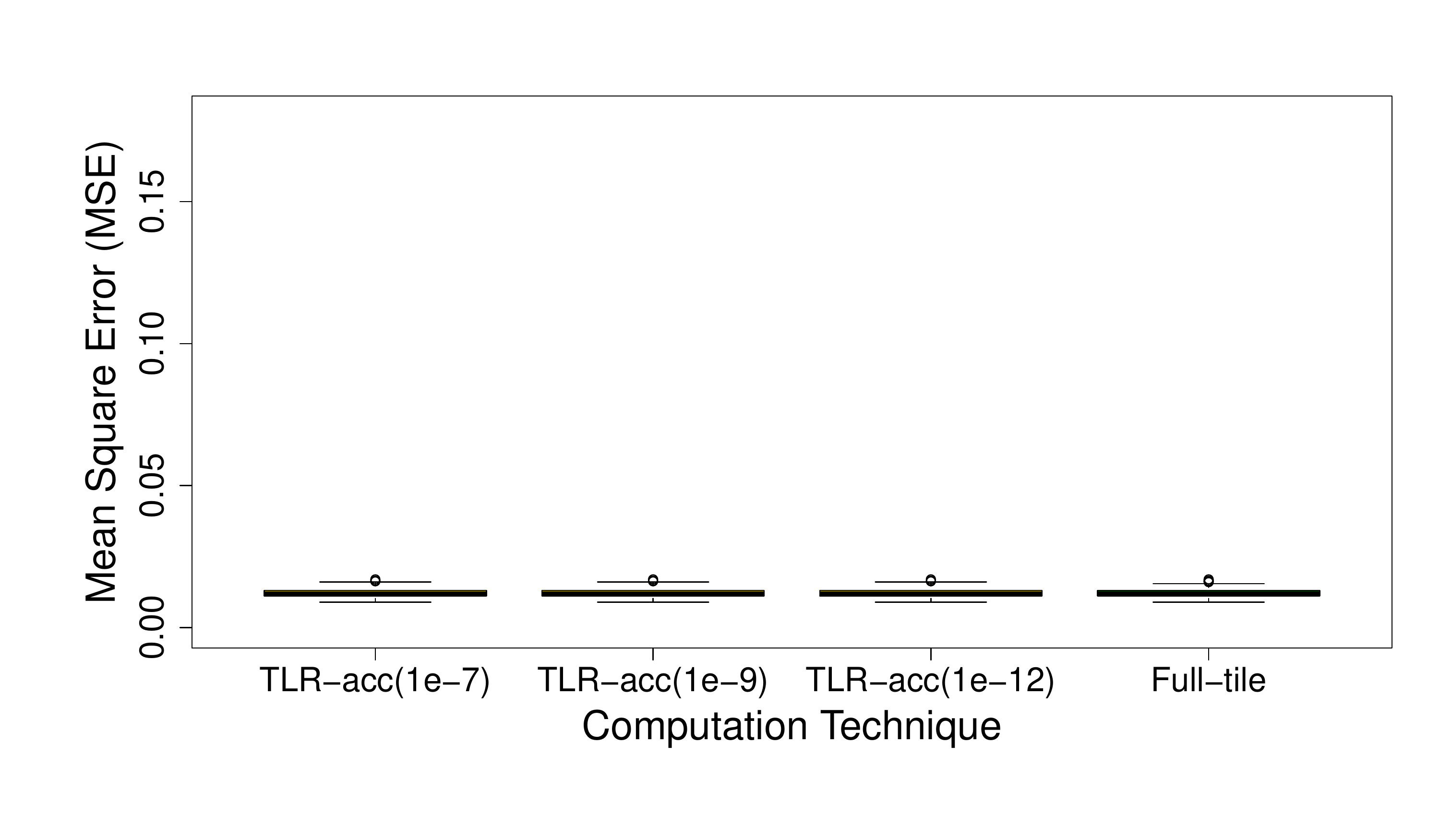}}

	\caption{Prediction mean square error (MSE) using synthetic datasets with three different parameter vector.}
	\label{fig:prediction-synthetic-boxplots}
\end{figure*}
The three boxplots are shown in Figure~\ref{fig:prediction-synthetic-boxplots}.
The TLR approximation variant perform well using the three accuracy thresholds (i.e., $10^{-7}$, $10^{-9}$, and $10^{-12}$) with
different parameter vector. This demonstrates the effectiveness of TLR approximation with
different data correlation degree in the prediction, even if the estimated parameter is not
as accurate as the \emph{full-tile} variant for some accuracy thresholds.

Another general observation for both TLR and \emph{full-tile} variants, prediction MSE becomes
lower in magnitude, if data is strongly correlated, as expected. For example, the average prediction MSE is $0.124$ in
the case of weak correlated data (i.e., $(1, 0.03, 0.5)$ ,  $0.036$ in the case of medium
correlated data (i.e., $(1, 0.1, 0.5)$), and $0.012$ in the case of strong correlated data (i.e., $(1, 0.3, 0.5)$).


\subsubsection{Real Datasets}
Qualitative assessment using real datasets is critical to ultimately assess the effectiveness of TLR approximation
for the MLE computations against the \emph{full-tile} variant. 
Here, we use the two different datasets, introduced in Section~\ref{sec:geospat}.

\begin{figure*}[!ht]
	\centering
		\subfigure[Soil moisture data ($8$ geographical regions).]{
			\label{fig:soil}
		\includegraphics[width=0.4\textwidth]{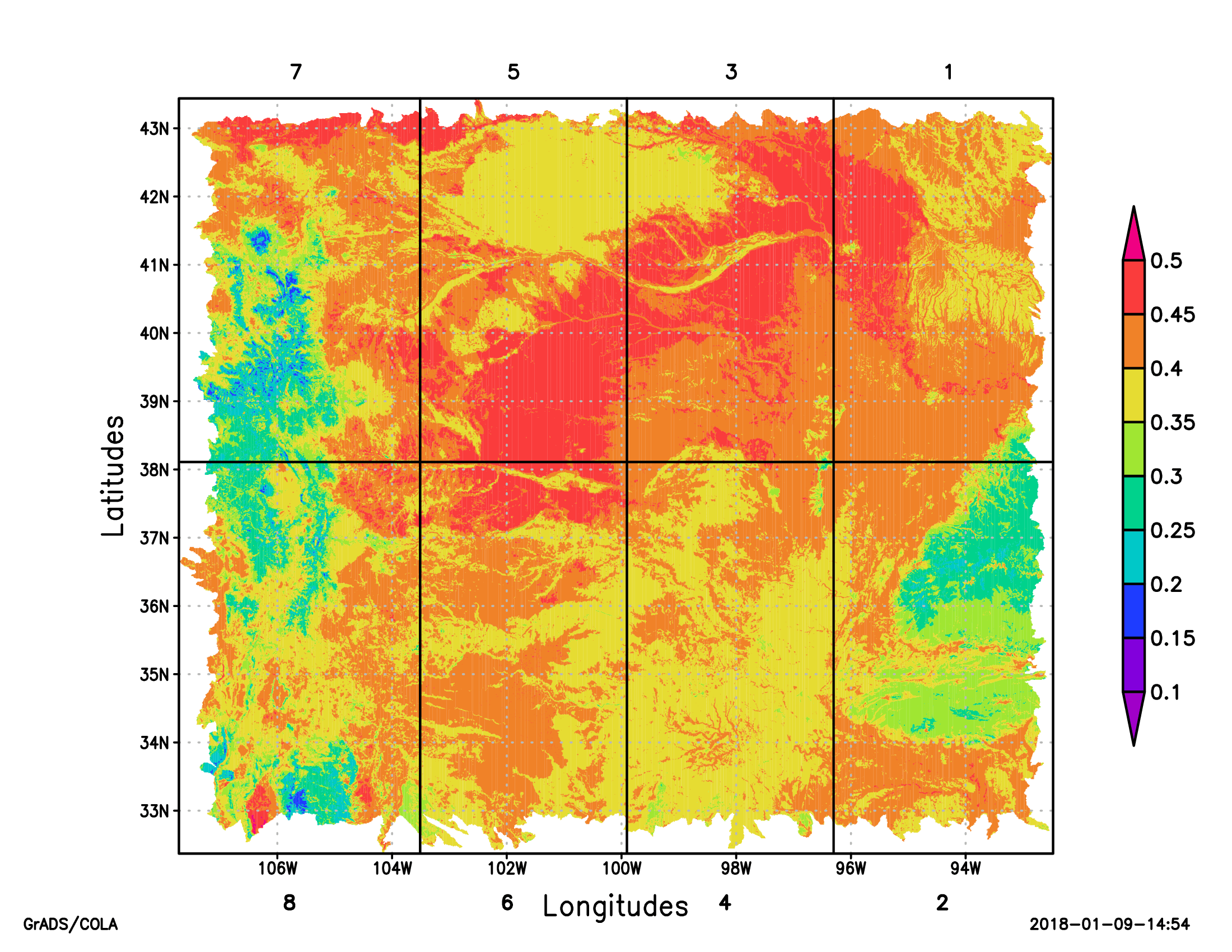}}
		\hspace{5mm}
		\subfigure[Wind speed data ($4$ geographical regions).]{
			\label{fig:wind}
		\includegraphics[width=0.4\textwidth]{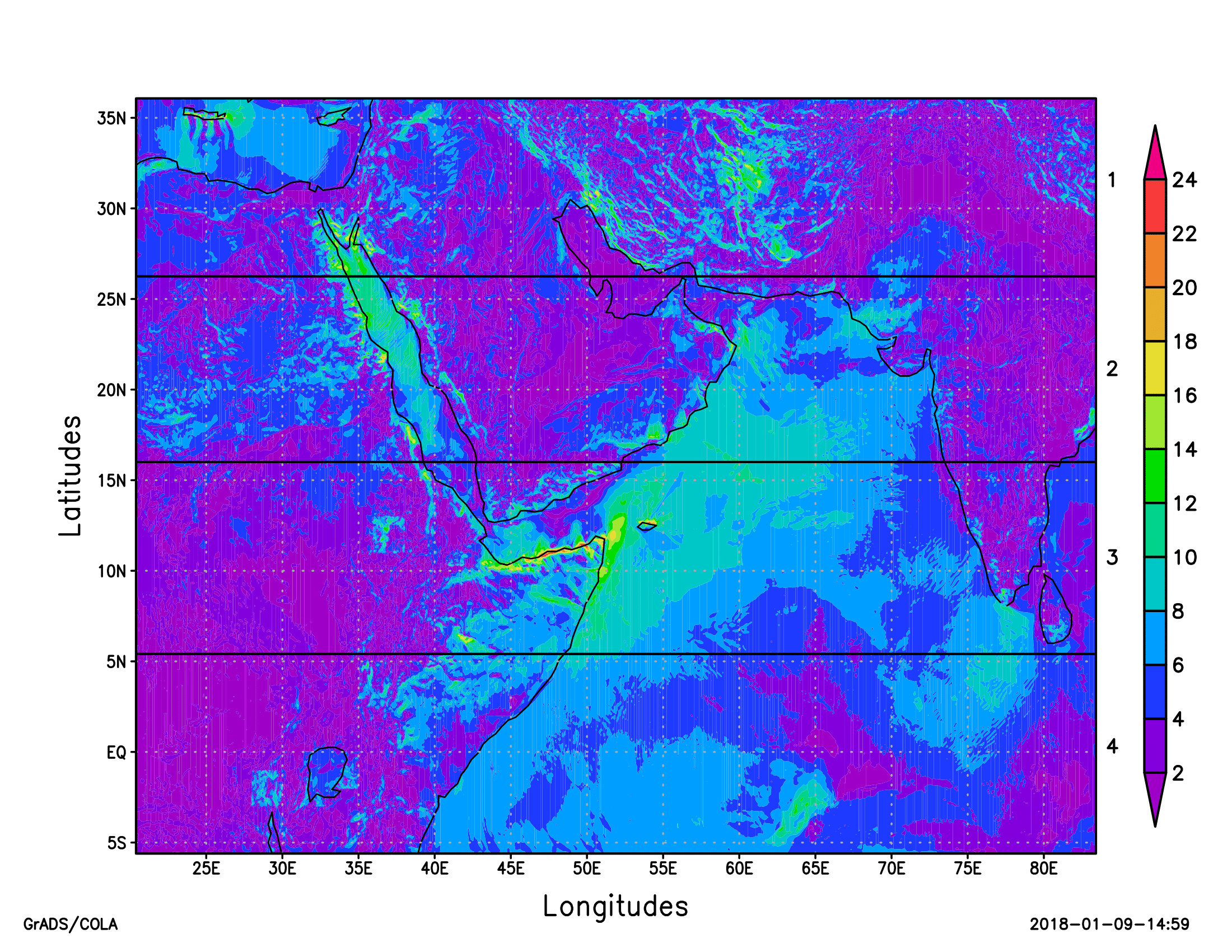}}
	\caption{Two examples of real geospatial datasets.}
	\label{fig:dataset}
\end{figure*}

Figure~\ref{fig:soil} shows the soil moisture dataset with $2$M locations. We divide
the data map into eight regions, from R$0$ to R$7$ to reduce the execution time of estimating the MLE
 operation especially in the case of exact computation. Furthermore, Figure~\ref{fig:wind} shows the
  wind speed dataset with $1$M locations. As the soil moisture dataset, we chose to divide
   the wind speed map into four regions from R$0$ to R$3$. 
   Generally, in both maps, each region contains about $250$K locations.

Tables~\ref{tab:regions8} and~\ref{tab:regions4} record the estimated parameters using TLR 
 approximation techniques with different accuracy thresholds as well as the reference one
 obtained with \emph{full-tile} variant. We report the estimated values
 to facilitate the reproducibility of this experiment. Using soil moisture datasets, we can use a TLR accuracy up to $10^{-12}$ which is
 still faster than \emph{full-tile}  technique while in wind speed dataset the highest accuracy used is $10^{-9}$ to maintain better performance
 compared to \emph{full-tile} variant.

Both tables show that highly correlated regions require high
TLR accuracy thresholds to reach the same parameter estimation quality as the \emph{full-tile} variant, 
e.g., the soil moisture data, R$7$ and R$8$, and 
the wind dataset, R$1$, R$2$, and R$3$. Moreover, the results show that the smoothness parameter is 
the easiest parameter to be estimated with any TLR accuracy thresholds, even 
in presence of highly correlated data.

\begin{table*}
\centering
	\setstretch{0.9}
\scriptsize
\caption{Estimation of the Mat\'{e}rn covariance parameters for $8$ geographical regions of the soil moisture dataset.}
\renewcommand{\arraystretch}{1.3}
\begin{tabu}{|c|[1pt]c|c|c|c|c|[1pt]c|c|c|c|c|[1pt]c|c|c|c|c|}
\hline
\scriptsize
& \multicolumn{15}{c|}{Mat\'{e}rn Covariance}  \\
R & \multicolumn{5}{c|[1pt]}{Variance ($\theta_1$)}  &  \multicolumn{5}{c|[1pt]}{Spatial Range ($\theta_2$)}  &  \multicolumn{5}{c|}{Smoothness ($\theta_3$)} \\  

& \multicolumn{4}{c|[0.5pt]}{TLR Accuracy}  &&  \multicolumn{4}{c|[0.5pt]}{TLR Accuracy} & &  \multicolumn{4}{c|[0.5pt]}{TLR Accuracy} & \multicolumn{1}{c|[0.5pt]}{} \\ 
& 

\scriptsize
   $10^{-5}$ & \scriptsize  $10^{-7}$ &   \scriptsize $10^{-9}$ &
     \scriptsize $10^{-12}$ &  \scriptsize \emph{Full-tile} &   
  \scriptsize $10^{-5}$ & \scriptsize  $10^{-7}$ &  \scriptsize $10^{-9}$ &
    \scriptsize  $10^{-12}$ &  \scriptsize \emph{Full-tile} &   
  \scriptsize $10^{-5}$ &  \scriptsize $10^{-7}$ &  \scriptsize $10^{-9}$ &   
  \scriptsize $10^{-12}$ & \scriptsize  \emph{Full-tile} 
 \\[-1pt] \tabucline[1pt]{1-19}
 \footnotesize
R1&
0.855&0.855&0.855&0.855&0.852&
6.039&6.034&6.034&6.033&5.994&
0.559&0.559&0.559&0.559&0.559\\ \hline
R2&
0.383&0.378&0.378&0.378&0.380&
10.457&10.307&10.307&10.307&10.434&
0.491&0.491&0.491&0.491&0.490\\ \hline
R3&
0.282&0.283&0.283&0.283&0.277&
11.037&11.064&11.066&11.066&10.878&
0.509&0.509&0.509&0.509&0.507\\ \hline
R4&
0.382&0.38&0.38&0.38&0.41&
7.105&7.042&7.042&7.042&7.77&
0.532&0.533&0.533&0.533&0.527 \\ \hline
R5&
0.832&0.837&0.837&0.837&0.836&
9.172&9.225&9.225&9.225&9.213&
0.497&0.497&0.497&0.497&0.496 \\ \hline
R6&
0.646&0.615&0.621&0.621&0.619&
10.886&10.21&10.317&10.317&10.323&
0.521&0.524&0.524&0.524&0.523\\ \hline
R7&
0.430&0.452&0.452&0.452&0.553&
14.101&15.057&15.075&15.075&19.203&
0.519&0.516&0.516&0.516&0.508\\ \hline
R8&
0.661&1.194&0.769&0.769&0.906&
18.603&37.315&22.168&22.168&27.861&
0.469&0.462&0.467&0.467&0.461\\ 
\hline

	\end{tabu}
	\label{tab:regions8}
\end{table*}

\begin{table*}
\centering
	\setstretch{0.9}
\scriptsize

\caption{Estimation of the Mat\'{e}rn covariance parameters for $4$ geographical regions of wind speed dataset.}
\renewcommand{\arraystretch}{1.3}
\begin{tabu}{|c|[1pt]c|c|c|c|[1pt]c|c|c|c|[1pt]c|c|c|c|}
\hline
& \multicolumn{12}{c|}{Mat\'{e}rn Covariance}  \\
R & \multicolumn{4}{c|[1pt]}{Variance ($\theta_1$)}  &  \multicolumn{4}{c|[1pt]}{Spatial Range ($\theta_2$)}  &  \multicolumn{4}{c|}{Smoothness ($\theta_3$)} \\

& \multicolumn{3}{c|[0.5pt]}{TLR Accuracy}  &&  \multicolumn{3}{c|[0.5pt]}{TLR Accuracy} &&  \multicolumn{3}{c|[0.5pt]}{ TLR Accuracy} & \multicolumn{1}{c|[0.5pt]}{} \\

& \scriptsize   $10^{-5}$ & \scriptsize  \scriptsize $10^{-7}$ & \scriptsize  $10^{-9}$  & \scriptsize    \emph{Full-tile} &\scriptsize  $10^{-5}$ & \scriptsize  $10^{-7}$ &  \scriptsize $10^{-9}$ &  \scriptsize   \emph{Full-tile} &  
\scriptsize $10^{-5}$ &  \scriptsize  $10^{-7}$ &   \scriptsize $10^{-9}$  & \scriptsize  \emph{Full-tile}
 \\[-1pt] \tabucline[1pt]{1-19}

R1	& 7.406		& 9.407	& 	12.247	& 	8.715	& 	29.576	& 	33.886		& 39.573	& 	32.083	& 	1.214& 1.196		& 1.175		& 	1.210\\ \hline
R2	& 11.920&	13.159		& 13.550		& 12.517	&26.011		& 28.083		& 28.707	& 	27.237		& 1.290 & 1.267		& 1.260		& 1.274 \\ \hline
R3		& 10.588	& 	10.944		& 11.232		& 10.819			& 18.423	& 	18.783	& 	19.114		& 18.634	&1.418	& 1.413		& 1.407			& 1.416\\ \hline
R4		& 12.408	& 	17.112		& 12.388	& 	12.270		& 17.264		& 17.112	& 	17.247		& 17.112	&  1.168 & 	1.170		& 1.168		& 1.170 \\ \hline

	\end{tabu}
\label{tab:regions4}
\end{table*}

\begin{figure*}[!ht]
  \centering
  \subfigure[Soil moisture data R1.]{
    \label{fig:theta1-soil8-region-0}
   \includegraphics[width=0.23\linewidth]{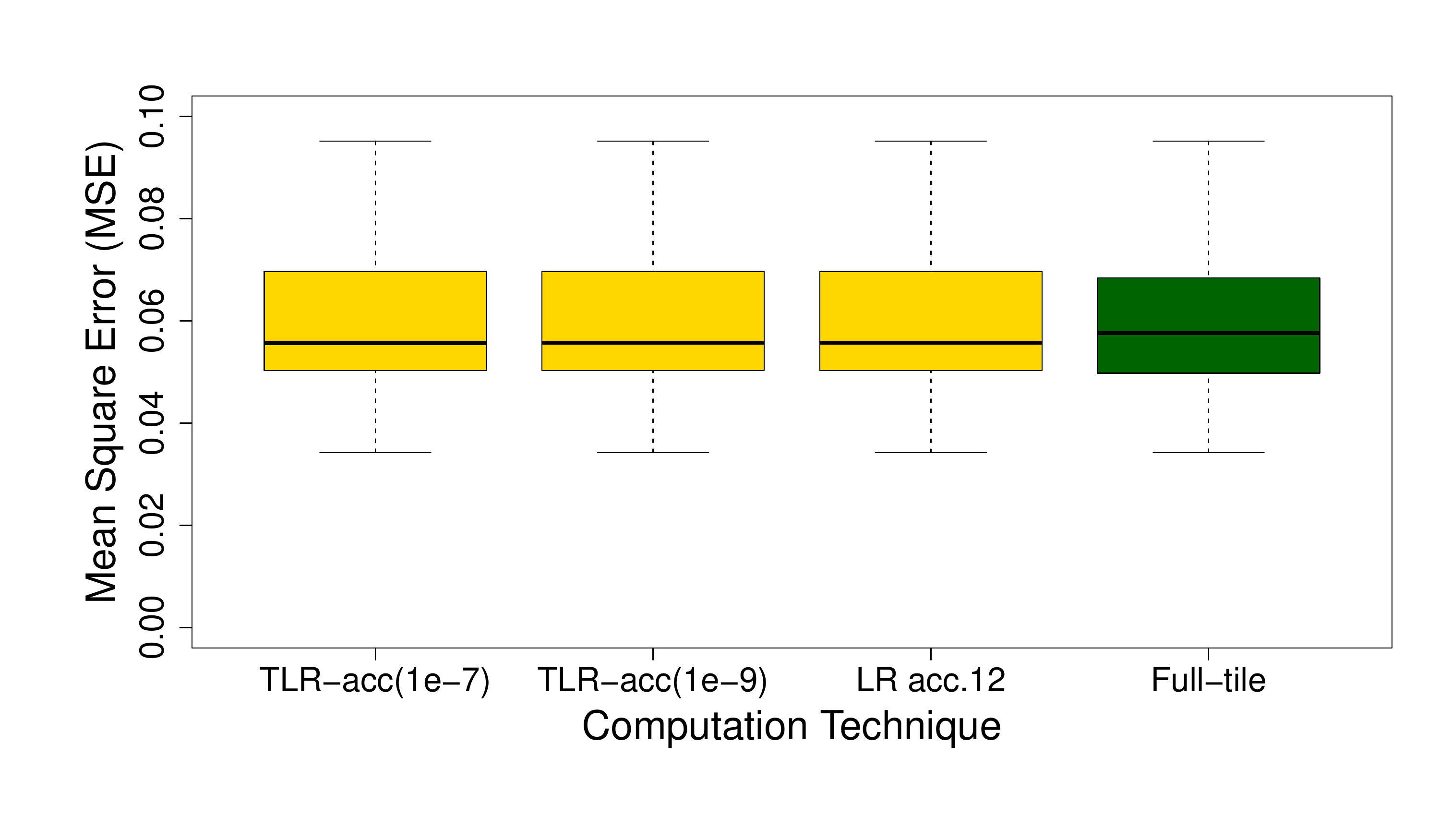}}
  \subfigure[Soil moisture data R3.]{
    \label{fig:theta2-soil8-region-3}
   \includegraphics[width=0.23\linewidth ]{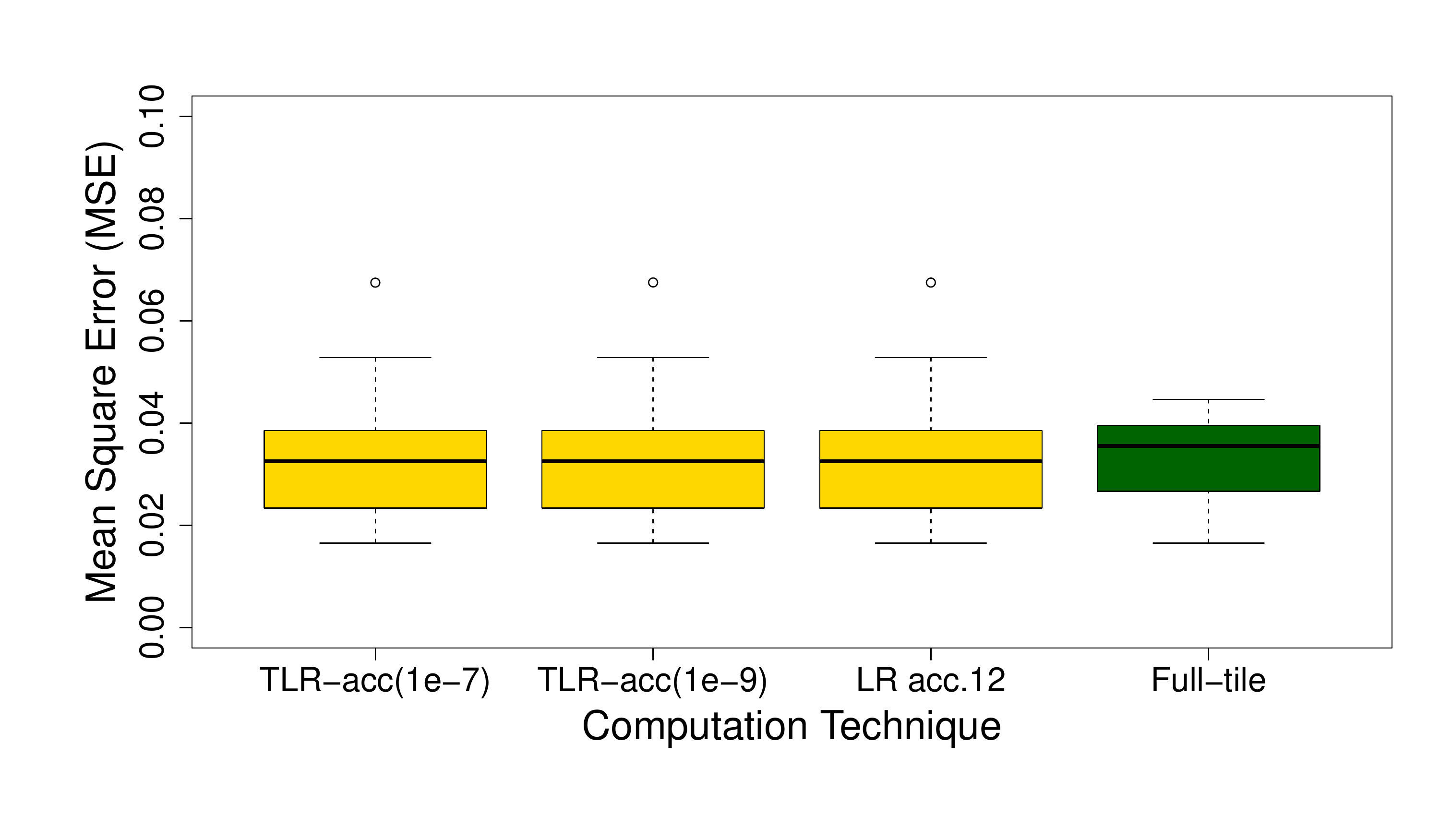}}
  \subfigure[Wind speed data R1.]{
    \label{fig:theta3-wind8-region-1}
   \includegraphics[width=0.23\linewidth ]{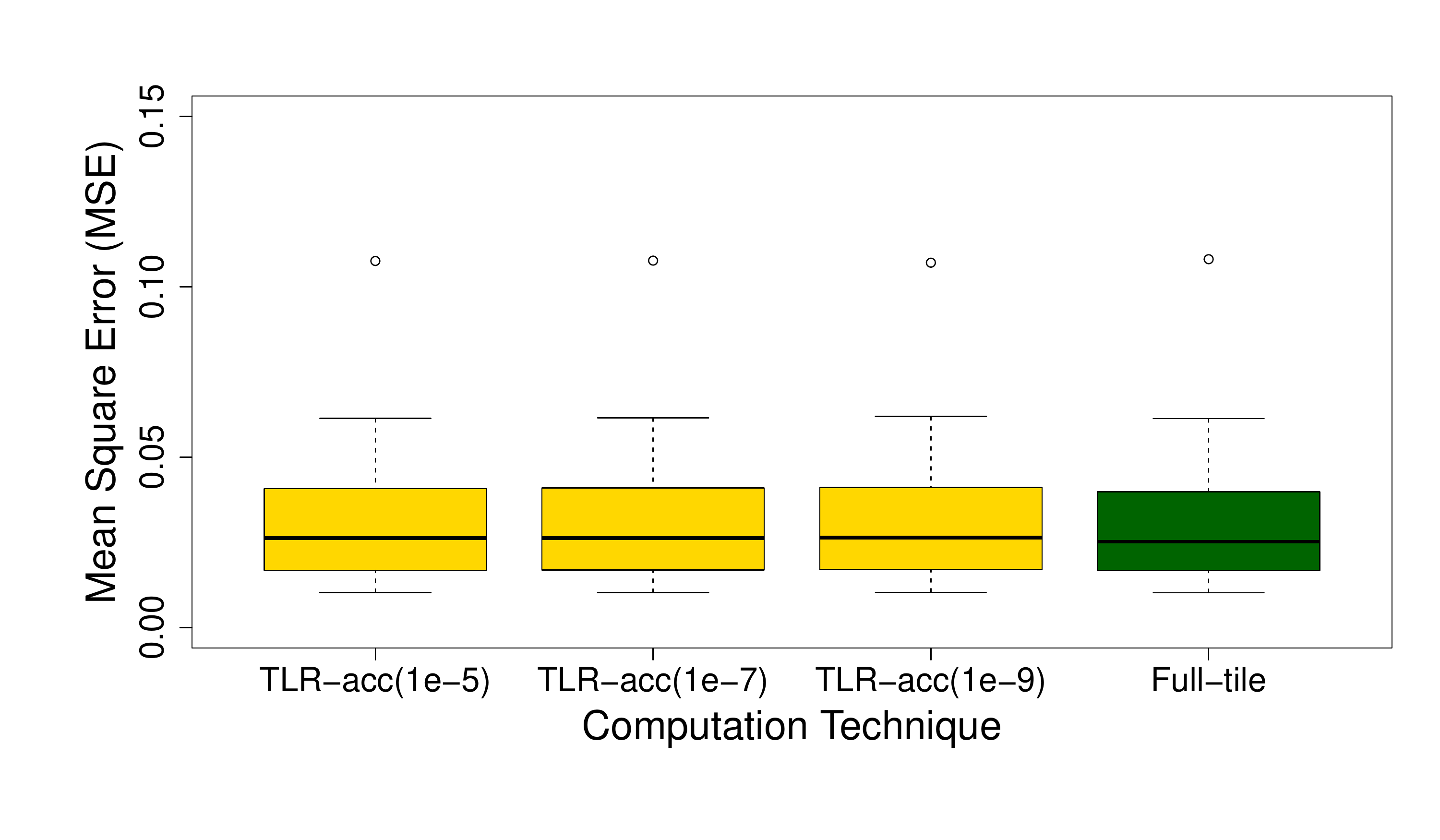}}
       \subfigure[Wind Speed data R4.]{
    \label{fig:theta1-wind8-region-3}
  \includegraphics[width=0.23\linewidth]{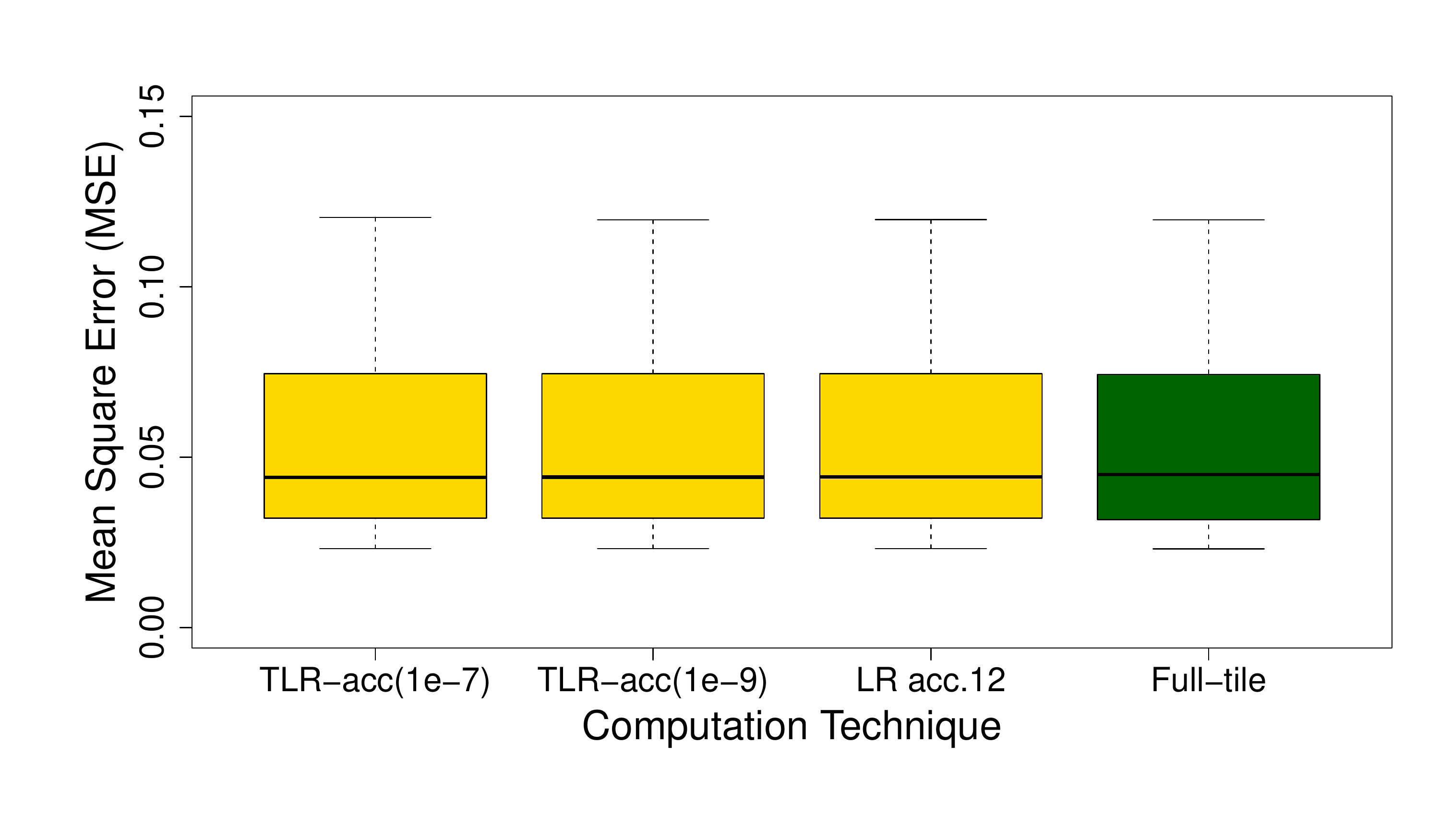}}
	\caption{ Prediction Mean Square Error (MSE) using Synthetic Datasets with three different parameter vector.}
	\label{fig:prediction-real-boxplots}
\end{figure*}

Moreover, we estimate the prediction MSE of $100$ missing values, which are randomly chosen from the same region. 
We select two regions from each dataset, i.e., R$1$ and R$4$ from the soil moisture data, and R$1$ and R$3$ from the wind speed data. 
We conduct this experiment $100$ times and Figure~\ref{fig:prediction-real-boxplots} 
shows the four boxplots with different computation techniques. The figure shows 
that TLR approximation technique for MLE provides a prediction 
MSE close to the \emph{full-tile} variant with 
different accuracy thresholds, even if the estimated parameters slightly differ, 
as shown in Tables~\ref{tab:regions8} and ~\ref{tab:regions4}.

\section{Conclusion}
\label{sec:summary}
This paper introduces the Tile Low-Rank approximation into the
open-source \texttt{ExaGeoStat} framework (\emph{TLR support to 
be released soon}) for effectively computing the 
Maximum Likelihood Estimation (MLE) on various parallel shared 
and distributed-memory systems, in the context of climate
and environmental applications. This permits to reduce 
the arithmetic complexity and memory footprint of MLE 
computations by exploiting the data sparsity structure of 
the Mat\'{e}rn covariance matrix of size up to $2$M.
The resulting TLR approximation for the MLE computation outperforms 
its full machine precision accuracy counterpart up to $13$X and $5$X
on synthetic and real datasets, respectively. A comprehensive qualitative assessment
of the accuracy of the statistical parameter estimation as well as the 
prediction (i.e., supervised learning) demonstrates the limited compromise
required to achieve high performance, while maintaining proper accuracy. 
%



\textbf{Acknowledgment.} 
We would like to thank 
Intel for support in the form of an Intel Parallel Computing Center award and Cray for support 
provided during the Center of Excellence award to the Extreme Computing
Research Center at KAUST.
This research made use of the resources of the KAUST Supercomputing Laboratory.

\IEEEtriggeratref{19}

\bibliographystyle{IEEEtran}
\bibliography{references}

\end{document}